\documentclass[10pt,twocolumn,twoside]{IEEEtran}
\usepackage{color,graphics,graphicx,amsmath,amsfonts,
amssymb,times,setspace,cite,fancyhdr,epstopdf}
\usepackage[active]{srcltx}

\usepackage[english]{babel} 
\usepackage{nomencl}
\usepackage{adjustbox}
\usepackage[dvipsnames]{xcolor}
\usepackage{enumerate}

\usepackage{array}
\usepackage{makecell}
\usepackage{booktabs}
\usepackage{multirow}
\usepackage{makecell}
\usepackage[utf8]{inputenc}
\usepackage{algpseudocode,algorithm,algorithmicx}
\makenomenclature
\usepackage{caption}
\usepackage{subcaption}

\def\BibTeX{{\rmfamily B\kern-.05em{\scshape i\kern-.025em b}\kern-.08em \TeX}}

\newcommand{\ben}{\begin{equation}}
\newcommand{\een}{\end{equation}}
\newcommand{\beq}{\begin{eqnarray}}
\newcommand{\eeq}{\end{eqnarray}}

\renewcommand\floatsep{0mm}
\renewcommand\textfloatsep{2mm}

\begin{document}

\title{Locational Scenario-based Pricing in a Bilateral Distribution Energy Market under Uncertainty}
\author{Hien Thanh Doan, Minsoo Kim, Keunju Song and Hongseok Kim}

\maketitle
\begin{abstract}

In recent years, there has been a significant focus on advancing the next generation of power systems. Despite these efforts, persistent challenges revolve around addressing the operational impact of uncertainty on predicted data, especially concerning economic dispatch and optimal power flow. To tackle these challenges, we introduce a stochastic day-ahead scheduling approach for a community. This method involves iterative improvements in economic dispatch and optimal power flow, aiming to minimize operational costs by incorporating quantile forecasting. Then, we present a real-time market and payment problem to handle optimization in real-time decision-making and payment calculation. We assess the effectiveness of our proposed method against benchmark results and conduct a test using data from 50 real households to demonstrate its practicality. Furthermore, we compare our method with existing studies in the field across two different seasons of the year. In the summer season, our method decreases optimality gap by 60\% compared to the baseline, and in the winter season, it reduces optimality gap by 67\%. Moreover, our proposed method mitigates the congestion of distribution network by 16.7\% within a day caused by uncertain energy, which is a crucial aspect for implementing energy markets in the real world.
\end{abstract}

%
\IEEEpeerreviewmaketitle

\section{Introduction} 
The strategic development of smart energy markets focused on decarbonization becomes essential when confronting the challenges presented by solar energy production variability. This challenge is underscored by the growing adoption of small-scale renewable energy sources, particularly photovoltaic (PV) systems. Despite the notable advantages offered by solar PV systems, such as reduced operational costs and environmental pollution mitigation, their reliance on solar irradiation introduces intricacies. The fluctuations in power generation not only compromise the power system's reliability but also introduce uncertainties in scheduling and an increased demand for reserve capacities. Consequently, this study focuses on designing smart energy markets to facilitate the seamless integration of renewable energy, optimize resource allocation, and mitigate the impact of uncertain energy, ultimately contributing to decarbonized and sustainable energy future.

In the past decade, solar PV systems have been the subject of extensive research across various domains. These investigations have delved into diverse areas, ranging from assessing the risks and benefits associated with the coordination of renewable energy sources (RESs) \cite{tan2021evaluation} to optimizing the planning of rooftop PV systems \cite{jung2021optimal}. Other studies have focused on minimizing costs for a large number of PV prosumers \cite{lopez2023peer}, analyzing energy planning and renewable energy plans \cite{icaza2022analysis}, optimizing the location planning of PV installations within transmission networks \cite{komiyama2019optimal}, and incorporating solar PV into transmission expansion planning \cite{somchit2023optimal,thongbouasy2023transmission}. Additionally, research has been conducted on microgrid energy management systems \cite{suresh2020microgrid} and flexible renewable energy planning \cite{kim2023flexible}.

Nevertheless, the primary challenge lies in addressing the uncertainties associated with the output of solar PV systems and household energy demand. To tackle these uncertainties, the widespread recommendation and application of energy storage systems (ESSs) have emerged \cite{kim2023stochastic}. ESSs have proven effective in mitigating the uncertainty levels of solar energy \cite{luo2020optimal}, exploring various types of ESSs, sizes, and locations while considering their impact on both active and reactive power \cite{lu2020data}. Specifically, research has delved into fair allocation methods for rolling horizon-based energy management systems \cite{erdincc2023rolling}, optimal dispatching strategies for uncertain demand-ESS scenarios \cite{zhou2023novel}, addressing minimum net load and operation costs \cite{huang2021investigation}, enhancing the hosting capacity of renewable energy system owners and system operations \cite{cao2021risk}, load frequency control for interconnected power systems \cite{gorbachev2023mpc}, and optimizing power flow control strategy of battery-PV for transportation systems \cite{ge2022combined}.

Indeed, the deployment of ESSs offers a potential solution for alleviating congestion on transmission lines under PV systems and demand. This has been explored from various angles, including the application of model predictive control (MPC) to minimize both grid congestion and ESS degradation \cite{nair2021analysis}, utilizing demand response (DR) programs to determine effective ESS strategies \cite{salehi2019optimal}, examining the impact of ESS on congestion with uncertain PV \cite{prajapati2021reliability}, investigating the coordination of PV-ESS to alleviate congestion \cite{wang2023benefits}, and optimal sizing for transmission congestion relief \cite{mazaheri2022data}.

However, existing studies \cite{nair2021analysis, salehi2019optimal, prajapati2021reliability, wang2023benefits, mazaheri2022data} could not overcome the uncertainties of PV generation and load fluctuation in running distribution energy markets. Specifically, in the current prosumer model, bilateral energy trade is not included uncetainties, which is a crucial component of energy markets \cite{prajapati2021reliability}. Furthermore, alternating current optimal power flow (AC-OPF) still poses challenges since it has not been comprehensively addressed, particularly when it comes to reactive power \cite{nair2021analysis}. There is also no assessment of how to impose prosumer (or equivalently called household hereafter) energy exchange costs in the existing model. Given the highlighted shortcomings in previous studies, the main contributions of this paper can be summarized as follows:
\begin{enumerate}
    \item We propose a bilateral scheduling framework for a ESS-based energy market aimed at minimizing the operational costs associated with the uncertainties of household PV generation and demand. This is considered under the perspectives of both virtualization (day-ahead) and realization (real-time and payment).
    \item Distribution energy market is introduced through a mathematical representation in the form of a mixed-integer conic programming. This model accounts for the uncertainties in community households by incorporating a deep-learning-based quantile forecasting model. Additionally, the integration of community owner and distribution system operator is modeled as an iterative stochastic optimization problem to alleviate congestion and minimize exchange losses costs, in contrast to the previous research that overlooked such costs.
    \item We introduce a novel pricing mechanism called integrated grid locational pricing (IGLP), which determines locational prices by leveraging distribution locational marginal pricing (DLMP) in distribution network. This pricing mechanism serves as a framework for day-ahead optimization, aiming to identify optimal values for real-time optimization. 
    \item We rigorously verify the effectiveness and superiority of the proposed methodology through the analysis of 50 real-household datasets and assessment of the impact of forecasting uncertainty on the energy market. The obtained results demonstrate that our method decreases the optimality gap by up to 67\% and effectively mitigates the congestion by 16.7\% in power lines within a day, unlike the high optimality gap and congestion observed in existing studies.
\end{enumerate}
\begin{figure*}
\centering
\includegraphics[width=0.7\textwidth]{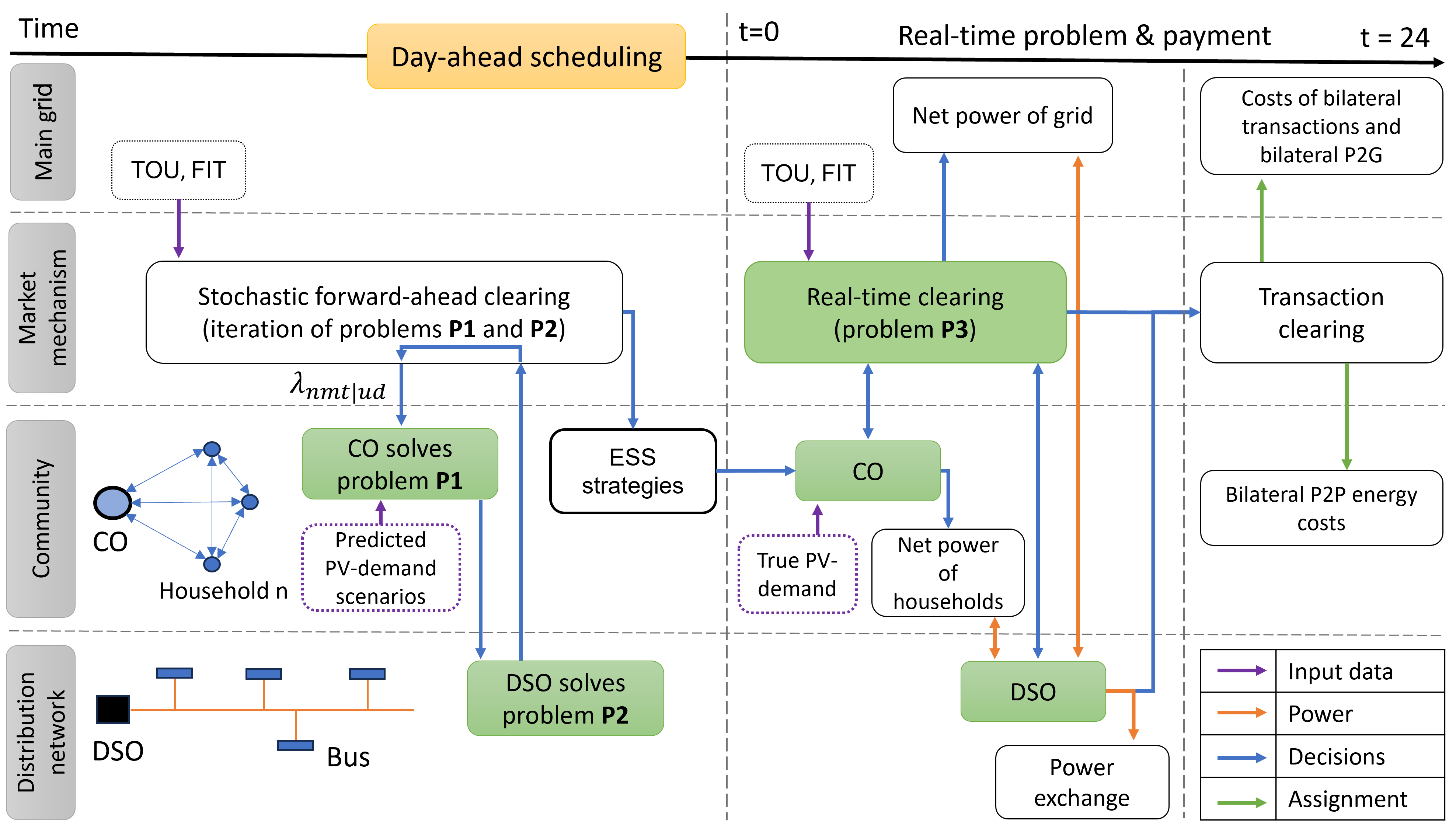}
\caption{A sequential market framework between community owner (CO) and distribution system operator (DSO), illustrates from left to right the day-ahead, real-time, and payment problem.}
\label{fig:overview}
\end{figure*}

The rest of this paper is organized as follows. Section~\ref{sec:system} describes the system model, and Section~\ref{sec:stochastic-optimization} proposes a method for addressing uncertainties in the day-ahead market. Following that, in Section~\ref{sec:Energy-market-clearing}, an alternative approach will be discussed for managing real-time market dynamics, alongside a discussion on the revision of the payment. Section~\ref{sec:performance} then presents the results and analyzes of the case studies, followed by the conclusion in Section~\ref{sec:conclusion}.
\section{System Model}\label{sec:system}
\subsection{Structure of energy market}
In this section, we delineate three main components of the proposed energy market: the community owner (CO), the distribution system operator (DSO), and the main grid, as illustrated in Fig.~\ref{fig:overview}. COs aim to optimize energy costs among households, manage communication exchanges with other elements and handle payment calculations. The DSO is responsible for monitoring and coordinating the reliability and stability of buses in the distribution network. The main grid serves as an auxiliary service for balancing the mismatch between generation and demand.

The participant in the market is a household. Let $\mathcal{N}$ represent a set of households indexed by $n=1,\ldots,|\mathcal{N}|$, and $\mathcal{T}$ denote a set of time slots in a day indexed by $t=1,\ldots,|\mathcal{T}|$. These households actively participate in the energy market with the goal of minimizing their energy costs over time slots $\mathcal{T}$. The energy market process can be divided into sequential decision-making problems within the framework: a day-ahead problem, a real-time problem, and a payment problem. The day-ahead problem is initiated in advance to estimate and address all potential outcomes based on uncertain PV generation and demand. The real-time market is activated when actual energy becomes available, and the payment process occurs later when the exchanged energy is recorded. 
\subsection{Energy community model}
Households are equipped with PV panels, electrical assets, ESS, and a smart meter, enabling them to participate in buying and selling energy, record the energy exchanged.
\subsubsection{Energy storage systems model}
Employed to store energy for later use, (\ref{eq:battery_SOC}) serves as the stored energy in battery ($E_{nt}$) of a household $n$ at time slot $t$. When the PV system supplies excess electricity, the household charges its battery $P_{nt}^{bc}$, leading to an increase in $E_{nt}$. Conversely, $E_{nt}$ decreases when the battery provides electricity, discharging the battery $P_{nt}^{bd}$ to meet demand or sell surplus energy to others. The efficiency of charging and discharging, denoted by $\eta$, is considered. Moreover, the charging $P_{nt}^{bc}$ and discharging $P_{nt}^{bd}$ of the battery are constrained by the maximum rate $P_n^b$ in (\ref{eq:battery_charge_limit}) and (\ref{eq:battery_discharge_limit}). As denoted in (\ref{eq:battery_SOC_limit}), $E_{nt}$ remains within the boundaries of $\underline{E}_n$ and $\overline{E}_n$, which are defined based on the specified usable battery capacity. The constraint (\ref{eq:battery_binary}) is enforced by binary variable $\psi_{nt}$ to prevent simultaneous charging and discharging. The battery constraints for household $n$ are summarized as follows:
\begin{subequations}
\label{eq:battery}
\begin{align}
E_{nt+1}= E_{nt}+\biggr[\eta P_{nt}^{bc} - P_{nt}^{bd}/\eta \biggr] \bigtriangleup_t, \label{eq:battery_SOC}
\end{align}
\begin{equation}
\underline {E}_n \leq E_{nt} \leq \overline {E}_n, \label{eq:battery_SOC_limit}
\end{equation}
\begin{equation}
0 \leq P_{nt}^{bc} \leq \psi_{nt} P_n^b, \label{eq:battery_charge_limit}
\end{equation}
\begin{equation}
0 \leq P_{nt}^{bd} \leq (1-\psi_{nt}) P_n^b, \label{eq:battery_discharge_limit}
\end{equation}
\begin{equation}
\psi_{nt} \in [0,1],  \label{eq:battery_binary}
\end{equation}
\end{subequations}
where $\bigtriangleup_t$ is a duration of each time slot.
\subsubsection{Peer-to-peer trades}
In the market, transactions involve peer-to-peer (P2P) exchanges, and the trade process can be decomposed into bilateral peer-to-peer setups \cite{frolke2022network}. This structure allows agents to directly negotiate and come to agreements regarding bilateral energy exchanges. The net trade $T_{nmt}$ of a household $n$ with its partners $m$ is expressed in (\ref{eq:net_trade_p2p}). This net trade consists of selling $s_{nmt} \geq 0$ ($n$ sells to $m$) and buying $b_{nmt} \geq 0$ ($n$ buys from $m$). The constraint (\ref{eq:bilateral_p2p}) ensures trade reciprocity through the relevant constraint between two households, $n$ and $m$.
\begin{subequations}
\begin{equation}
    T_{nmt} = s_{nmt} - b_{nmt}, \label{eq:net_trade_p2p}
\end{equation}
\begin{equation}
    s_{nmt} = b_{mnt}. \label{eq:bilateral_p2p}
\end{equation}
\end{subequations}
\subsubsection{Energy balance model}
Every household has a main grid connection, peer-to-grid, which allows them to withdraw or inject electricity. The active power net injection of household $n$ at time slot $t$ into the power system is designated as $P_{nt}^P$, with positive values representing injections. (\ref{eq:house_net_active_power}) shows how $P_{nt}^P$ is calculated. The generation energy that a household $n$ has at time slot $t$ includes PV generation $G_{nt}^{pv}$, energy discharge from battery $P_{nt}^{bd}$. Meanwhile, the energy demand comprises the appliance demand $L_{nt}^{de}$, energy charge to battery $P_{nt}^{bc}$. Households can sell energy to grid $s_{n0t}$ (the grid is indexed by $0$ hereafter) and to other houses $s_{nmt}$, $m\in \omega_n$, where $\omega_n$ represents a set of households engaged in trading with household $n$, as indicated in (\ref{eq:net_injection_balance}). Additionally, households can purchase energy from the grid $b_{n0t}$ and from other houses $b_{nmt}$. Besides, reactive power $P_{nt}^Q$ can be calculated by (\ref{eq:house_net_reactive_power}):
\begin{subequations}
\label{eq:energy_balance}
\begin{equation}
    P_{nt}^P = G_{nt}^{pv} + P_{nt}^{bd} - L_{nt}^{de} - P_{nt}^{bc}, \label{eq:house_net_active_power}
\end{equation}
\begin{equation}
    P^P_{nt} = s_{n0t} + \sum_{m\in \omega_n} s_{nmt} - b_{n0t} - \sum_{m\in \omega_n} b_{nmt},
    \label{eq:net_injection_balance}
\end{equation}
\begin{equation}
    P_{nt}^{Q} = - L_{nt}^{Q}. \label{eq:house_net_reactive_power}
\end{equation}
\end{subequations}
\subsection{Distributed network model}
Consider a radial low-voltage (LV) distribution grid owned by DSO and denoted as $\mathcal{G}=(\mathcal{I},\mathcal{L})$, represented by an undirected connected graph, where $\mathcal{I}$ is the set of nodes $i=0,1,\ldots,|\mathcal{I}|$, and $\mathcal{L}$ is a set of connecting lines $\ell=1,\ldots,|\mathcal{L}|$. The slack bus, indexed as 0, is the root node. Each node $i$ has a parent node $\mathcal{A}_i$, and a set of children $\mathcal{C}_{i}$, where each child is indexed $j=1,\ldots,|\mathcal{C}_i|$ in a radial distribution network. The current flow from the parent $\mathcal{A}_i$ to node $i$ is denoted as $I_i$. The active and reactive powers of generation and demand are represented by $g_i^P$, $g_i^Q$, $d_i^P$, and $d_i^Q$, respectively. The active and reactive power flows on line $i$ are $f_i^P$ and $f_i^Q$. Note that, in a radial network, line $i$ is uniquely determined as the line from $\mathcal{A}_i$ to node $i$, so node index $i$ can be also used for line index without ambiguity. Each node has the squared voltage and current denoted as $v_i$ and $l_i$, respectively. The minimum and maximum squared voltage limits are $v_{i}^{\min}$ and $v_{i}^{\max}$. Each line is characterized by its resistance $r_i$ and reactance $x_i$, with maximum capacity $S^{\max}_{i}$. The conductance and susceptance of bus $i$ are $G_i$ and $B_i$, respectively. In the context of AC-OPF analysis for the radial network, these quantities are interconnected using the DistFlow approach, as described by the following constraints associated the corresponding Lagrange multipliers $(\gamma_{\mathcal{A}_i},\mu_i,\eta_{i}^{+},\eta_{i}^{-})$ for all $i$ in $\mathcal{I}$ \cite{kim2019p2p}. 
\begin{subequations}
\label{eq:ac-opf}
\begin{align}
(\gamma_{\mathcal{A}_i}) &: f_{it}^P + g_{it}^P-\sum_{j\in \mathcal{C}_i}(f_{jt}^P-r_{j}l_{jt}) + G_{i}v_{i} = d_{it}^P, \label{eq:ac-opf_active_flow}\\
(\mu_i) &: f_{it}^Q + g_{it}^Q-\sum_{j\in \mathcal{C}_i}(f_{jt}^Q-x_{j}l_{jt}) + B_{i}v_{i} = d_{it}^Q, \label{eq:ac-opf_reactive_flow}\\
(\eta_{i}^{+}) &: (f_{it}^P)^2+(f_{it}^Q)^2 \leq (S^{\max}_{i})^2, \label{eq:ac-opf_forward_limit} \\
(\eta_{i}^{-}) &: (f_{it}^P - l_{it}r_{i})^2+(f_{it}^Q - l_{it}x_{i})^2 \leq (S^{\max}_{i})^2, \label{eq:ac-opf_backward_limit}
\end{align}
\begin{align}
{}& v_{it}+2(r_{i}f_{it}^P+x_{i}f_{it}^Q)+l_{it}(r_{i}^2+x_{i}^2) = v_{\mathcal{A}_it}, \label{eq:ac-opf_voltage_balance}\\
{}& ||2f_{it}^P,2f_{it}^Q,v_{it}-l_{it}||_2 \leq v_{it}+l_{it}, \label{eq:ac-opf_SOC-constraint}
\end{align}
where the limit of the squared voltage $v_{i}$ of each bus $i$ can be written as,
\begin{align}
{}&\underline v_{i} \leq v_{it} \leq \overline v_{i}. \label{eq:ac-opf_vol_lim}
\end{align}
\end{subequations}
\subsection{Incorporation energy model}
In order to map the household owned by the CO and the bus in the power distribution system, we designate $\Omega_i$ as a set of the households $n$ that are connected to bus $i$. Then, the net injection active power of bus $i$ should be the total net injection of households $P_{nt}^P$ connected to bus $i$,  
\begin{equation}
    g_{it}^P - d_{it}^P = \sum_{n\in \Omega_{i}} P_{nt}^P \triangleq h_{it}^P,
\end{equation}
which is denoted by $h_{it}^P$. Similarly, the net injection reactive power in (\ref{eq:ac-opf_reactive_flow}) should be the total net injection of households $P_{nt}^Q$ located at bus $i$, denoted by $h_{it}^Q$. Then,
\begin{equation}
    - d_{it}^Q = \sum_{n\in\Omega_{i}} P_{nt}^Q \triangleq h_{it}^Q,
\end{equation}
Then, (\ref{eq:ac-opf_active_flow}) and (\ref{eq:ac-opf_reactive_flow}) can be rewritten by  
\begin{subequations}
\begin{align}
(\gamma_{\mathcal{A}_i}) &: f_{it}^P + h_{it}^P -\sum_{j\in \mathcal{C}_i}(f_{jt}^P-r_{j}l_{jt}) + G_{i}v_{it} &=& 0, \label{eq:ac-opf_net_active_flow}\\
(\mu_i) &: f_{it}^Q + h_{it}^Q -\sum_{j\in \mathcal{C}_i}(f_{jt}^Q-x_{j}l_{jt}) + B_{i}v_{it} &=& 0. \label{eq:ac-opf_net_reactive_flow}
\end{align}
\end{subequations}

Derived from the above equations, the energy market can be characterized by mixed-integer conic programming (MICP). In this framework, households try to minimize their energy costs through the optimization of self-consumption and participation in bilateral trade with others households in a mixed-integer linear problem (MILP). As a result, they decrease their dependency on the main grid, leading to significant cost savings while considering the optimal power flow in conic problem (CP). The formulation of this problem \textit{without} uncertainty is outlined below.

\vspace{2mm}
\noindent \textbf{P0: Deterministic energy cost minimization}
\begin{align} 
\text{minimize ~} {}& \sum_{n\in\mathcal{N}} \sum_{t\in\mathcal{T}} (\overline \lambda_{t} b_{n0t} - \underline \lambda_{t} s_{n0t} + \overline {\lambda}_{t}{g}_{0t}^{P}) \bigtriangleup_t \label{eq:SWM} \hspace{10000pt minus 1fil} \\
\text{subject to } {}& \text{ESS: (\ref{eq:battery_SOC})--(\ref{eq:battery_binary})}, \nonumber \\ 
{}& \text{P2P: (\ref{eq:net_trade_p2p})--(\ref{eq:bilateral_p2p}), (\ref{eq:house_net_active_power})--(\ref{eq:house_net_reactive_power})}, \nonumber \\ 
{}& \text{NET: (\ref{eq:ac-opf_forward_limit})--(\ref{eq:ac-opf_vol_lim}), (\ref{eq:ac-opf_net_active_flow})--(\ref{eq:ac-opf_net_reactive_flow})}, \nonumber \\
\text{variables ~} {}&  \{b_{nmt}, s_{nmt}, T_{nmt},P^P_{nt}, P_{nt}^{bd}, P_{nt}^{bc}, {g}_{0t}^{P}, \hspace{10000pt minus 1fil} \nonumber \\ 
{}& \psi_{nt}, b_{n0t}, s_{n0t}, v_{it}, l_{it}, h_{it}^P, h_{it}^Q, \nonumber \\
{}& n\in \mathcal{N}, i\in \mathcal{I}, n \in \Omega_i, m \in \omega_n, t\in \mathcal{T} \nonumber \}. 
\end{align}
Here, $\overline{\lambda}_{t}$ denotes the energy purchase price from the slack bus (main grid) located at node $0$, and $\underline{\lambda}_{t}$ represents the energy selling price to the grid. ${g}_{0t}^{P}$ represents the quantity of power losses imposed. The slack bus operates as an auxiliary service provider in the energy market. It sells energy to designated nodes when overall demand exceeds local generation and purchases energy when total generation surpasses local demand. The slack bus also provides reactive power to ensure power system reliability and stability.

\section{Stochastic optimization and day-ahead market}\label{sec:stochastic-optimization}
Centralized energy management in \textbf{P0} faces three main issues. Firstly, there is a privacy concern as the CO and the DSO have to share their private information \cite{tsaousoglou2022market}. Secondly, the exponential complexity of MICP $O(2^D)$ depends on number of decision variables $D$, significantly impacting the scalability and feasibility of the energy market \cite{liu2019distributed}. Thirdly, neglecting uncertainty may lead to obtaining sub-optimal solutions.

As a result, we divide Problem \textbf{P0} into two sub problems as shown in Fig.~\ref{fig:overview}. The first Problem \textbf{P1} in community level focuses on stochastic economic dispatch with the primary objective of maximizing social welfare while minimizing generation costs but ignores the distribution network constraints of $\text{NET: (\ref{eq:ac-opf_forward_limit})--(\ref{eq:ac-opf_vol_lim})}$ and $\text{(\ref{eq:ac-opf_net_active_flow})--(\ref{eq:ac-opf_net_reactive_flow})}$. The second problem \textbf{P2} in distribution network level exclusively deals with transmission cost minimization with the overarching goal of reducing transmission expenses in distributed manner. Then, \textbf{P1} and \textbf{P2} are connected by introducing a coupling variable named linear bilateral energy price under uncertainty. It streamlines the solution process for both problems as illustrated in Fig.~\ref{fig:overview} and seeks an optimal solution for the overarching problem. Moreover, iterative problem-solving is enabled, where the outcomes of the two problems are used to fine-tune the original problem solution by coupling variables.

\subsection{ Deep learning based scenario generation} \label{subsec:quantile_model}
Before describing stochastic P2P optimization, we introduce how to address the uncertainty associated with PV generation and demand. In doing this we utilize a deep learning based quantile forecasting model to generate scenarios for stochastic optimization. Specifically, we employ the model developed in our prior work of \cite{kim2023stochastic}. Unlike point-forecasting models that use mean squared error or mean absolute error as a loss function, the employed quantile forecasting model uses a quantile loss function which is defined as follows \cite{kim2023stochastic}:
\begin{equation}
\begin{aligned}
l_q({\hat{y}_q^t, y^t})
=\begin{cases}
(1-q)|\hat{y}_q^t-y^t|, & \text{if $y^t \leq \hat{y}_q^t$} \\
q|\hat{y}_q^t-y^t|, & \text{otherwise}
\end{cases}
\end{aligned}
\label{eq:loss_q}
\end{equation}
where $\hat{y}_q^t$ is the forecasted results at quantile $q\in(0,1)$ and $y^t$ is the corresponding labels in time slot $t$, respectively. Since the model with (\ref{eq:loss_q}) can only forecast a single quantile of a single time slot $t$, we adapt the loss function to generate multiple quantiles of scenarios with $T$ time slots as follows:
\begin{equation}
    L_q = \frac{1}{|\mathcal{Q}|} \sum_{t=1}^T\sum_{q\in \mathcal{Q}}l_q({\hat{y}_q^t, y^t}),
\label{eq:quantile_loss}
\end{equation}
where $\mathcal{Q}$ is a set of desired quantiles and $|Q|$ denotes the cardinality of $\mathcal{Q}$. Then, we use a set of the scenarios as $\mathcal{Y} = \{y_q | y_q = \{\hat{y}_q^t\}_{t=1}^T\in\mathbb{R}^T, \forall q\in\mathcal{Q}\}$ and select \textit{equally} spaced quantiles, i.e., $\mathcal{Q} = \{\delta, 2\delta, ..., 1 - \delta\}$ where $\delta = 1/(|\mathcal{Q}| + 1)$. In \cite{kim2023stochastic}, we proved that all the scenarios generated from this model have the same probability and serve as input for the stochastic P2P problem in the next subsection. Note that this model is general in the sense that $\hat{y}_q^t$ can be obtained by any deep learning based forecasting model, and we adopt TransLSTM, developed in our prior work \cite{song2022dttrans}.

\subsection{Stochastic P2P optimization problem} \label{subsec:PG_model}
The stochastic economic dispatch of the energy market can be formulated as an MILP problem as will be described in (\ref{eq:sto_CO_problem}). The main goal is to minimize the overall cost under the uncertainty induced by a set of PV scenarios denoted by $\mathcal{Y}^{pv}$ and a set of demand scenarios denoted by $\mathcal{Y}^{de}$. We exploit quantile scenarios $u\in\mathcal{Y}^{pv}$ for PV, and quantile scenarios $d\in\mathcal{Y}^{de}$ for demand. Note that \textbf{P1} mainly differs from \textbf{P0} in terms of addressing uncertainty, and scenarios are denoted by subscripts; for example, $b_{n0t|ud}$ denotes $b_{n0t}$ given the scenarios of $\mathcal{Y}^{pv}$ and $d\in\mathcal{Y}^{de}$. The virtue of employing quantile scenarios are such that all scenarios are equally likely, and thus the stochastic economic dispatch energy scheduling can be simply formulated as follows.

\vspace{2mm}
\noindent \textbf{P1: Stochastic economic dispatch of CO}
\begin{subequations}   
\label{eq:sto_CO_problem}
\begin{flalign}
\text{minimize } \sum_{n\in\mathcal{N}} \sum_{t\in\mathcal{T}} \sum_{u\in \mathcal{Y}^{pv}} \sum_{d\in \mathcal{Y}^{de}} \frac{1}{|\mathcal{Y}^{pv}||\mathcal{Y}^{de}|} \biggr(\overline \lambda_{t} b_{n0t|ud} \nonumber \\ - \underline \lambda_{t}s_{n0t|ud}
-\sum_{m\in \omega_n}\lambda_{nmt|ud}T_{nmt|ud} \biggr) \bigtriangleup_t 
\label{eq:sto_market_obj_fn}
\end{flalign}
\begin{equation}
\text{subject to~~~~~~~~~~~~~~~~~~~~~~~(\ref{eq:battery_SOC})-(\ref{eq:battery_binary}), } \hspace{10000pt minus 1fil} 
\end{equation}
\begin{equation}
    P_{nt|ud}^{P} = \hat{G}_{nt|u}^{pv} + P_{nt}^{bd} - \hat{L}_{nt|d} - P_{nt}^{bc}, \label{eq:sto_house_net_active_power}
\end{equation}
\begin{equation}
    P_{nt|ud}^{P} = s_{n0t|ud} + \sum_{m\in \omega_n} s_{nmt|ud} - b_{n0t|ud} - \sum_{m\in \omega_n} b_{nmt|ud}, \label{eq:sto_net_injection_active_power}
\end{equation}
\begin{equation}
    P_{nt|ud}^{Q} = - {L}_{nt|d}^{Q}, \label{eq:sto_house_net_reactive_power}
\end{equation}
\begin{equation}
    s_{nmt|ud}=b_{mnt|ud}, \label{eq:sto_bilateral_p2p}
\end{equation}
\begin{equation}
    T_{nmt|ud} = s_{nmt|ud} - b_{nmt|ud}, \label{eq:sto_net_trade_p2p}
\end{equation}
\begin{align}
\text{variables } {}&  \{b_{n0t|ud}, s_{n0t|ud}, T_{nmt|ud}, P_{nt|ud}^{P}, \hspace{10000pt minus 1fil} \nonumber \\ 
{}& P_{nt}^{bd}, P_{nt}^{bc},\psi_{nt}, s_{nmt|ud}, b_{nmt|ud}, \nonumber \\
{}& n\in \mathcal{N}, m\in \omega_{n}, u\in \mathcal{Y}^{pv}, d\in \mathcal{Y}^{de}, t\in \mathcal{T} \nonumber \},
\end{align}
\end{subequations}
where $\hat{G}^{pv}_{nt|u}$ denotes predicted PV generation scenario data, and $\hat{L}_{nt|d}$ denotes the predicted demand. The bilateral P2P price $\lambda_{nmt|ud}$ in (\ref{eq:sto_market_obj_fn}) depends on whether $n$ is a buyer or a seller and it serves as a \textit{coupling} variable between the CO and the DSO, and will be discussed in detail in \textbf{P2}.
\subsection{Distributed energy exchange costs problem} \label{subsec:DEXCP}
Since \textbf{P1} does not consider the distribution network constraints, the solution to the Problem \textbf{P1} may cause overload in distribution networks, resulting in the constraint violation in $\text{(\ref{eq:ac-opf_forward_limit})--(\ref{eq:ac-opf_vol_lim})}$ and $\text{(\ref{eq:ac-opf_net_active_flow})--(\ref{eq:ac-opf_net_reactive_flow})}$. To address this challenge, we introduce the concept of slack variables as a means of mitigating potential violations. The overarching objective in (\ref{eq:AC-OPF_scenario-objective}) of \textbf{P2} for the network's owner is twofold: to minimize both the energy cost and the additional cost associated with any constraint violations. To formalize this approach, the power flows and the distribution network constraints of (\ref{eq:ac-opf_forward_limit})--(\ref{eq:ac-opf_vol_lim}) and (\ref{eq:ac-opf_net_active_flow})--(\ref{eq:ac-opf_net_reactive_flow}) are reformulated to (\ref{eq:market-objective_active_flow})--(\ref{eq:market-objective_reactive_lim}) given scenarios. In doing this we newly introduce two slack variables $\xi_{it|ud}$ and $\varphi_{it|ud}$ in order to consider constraints violation of power flow and voltage, respectively.

\vspace{2mm}
\noindent \textbf{P2: AC-OPF of DSO}
\begin{subequations}
\begin{flalign}
\label{eq:AC-OPF_scenario-objective}
\text{minimize~~}  \sum_{t\in\mathcal{T}} \sum_{u\in \mathcal{Y}^{pv}} \sum_{d\in \mathcal{Y}^{de}} (\overline \lambda_t{g}_{0t|ud}^{P} + \varpi \xi_{it|ud} + \tau \varphi_{it|ud} ) \bigtriangleup_t
\end{flalign}
\begin{equation}
\hfilneg \text{subject to } \hspace{10000pt minus 1fil} \nonumber
\end{equation}
\begin{equation}
(\gamma_{\mathcal{A}_{i}t|ud}) : f_{it|ud}^P + h_{it|ud}^P -\sum_{j\in \mathcal{C}_i}(f_{jt|ud}^P-r_{j}l_{jt|ud}) + G_{i}v_{it|ud} = 0, \label{eq:market-objective_active_flow}    
\end{equation}
\begin{equation}
(\mu_{it|ud}) : f_{it|ud}^Q + h_{it|ud}^Q -\sum_{j\in \mathcal{C}_i}(f_{jt|ud}^Q-x_{j}l_{jt|ud}) + B_{i}v_{it|ud} = 0, \label{eq:market-objective_reactive_flow}
\end{equation}
\begin{flalign}
(\eta_{it|ud}^{+}) &: (f_{it|ud}^P)^2+(f_{it|ud}^Q)^2 \leq (S^{\max}_{i})^2 + \xi_{it|ud},& \label{eq:market-objective_forward_limit}    
\end{flalign}
\begin{flalign}
(\eta_{it|ud}^{-}) : (f_{it|ud}^P - l_{it|ud}r_{i})^2+(f_{it|ud}^Q - l_{it|ud}x_{i})^2 \nonumber \\  \leq (S^{\max}_{i})^2 + \xi_{it|ud}, \label{eq:market-objective_backward_limit}
\end{flalign}
\begin{equation}
v_{it|ud}+2(r_{i}f_{it|ud}^P+x_{i}f_{it|ud}^Q)+l_{it|ud}(r_{i}^2+x_{i}^2) = v_{\mathcal{A}_{i}t|ud}, \label{eq:market-objective_vol_inequal}
\end{equation}
\begin{equation}
 ||2f_{it|ud}^P,2f_{it|ud}^Q,v_{it|ud}-l_{it|ud}||_2 \leq v_{it|ud}+l_{it|ud},  \label{eq:market-objective_SOC-constraint}
\end{equation}
\begin{equation}
    0 \leq \underline v_{i} - \varphi_{it|ud} \leq v_{it|ud} \leq \overline v_{i} + \varphi_{it|ud}, \label{eq:market-objective_vol_lim}
\end{equation}
\begin{equation}
    \underline {g}_0^{P} \leq {g}_{0t|ud}^{P} \leq \overline {g}_0^{P}, \label{eq:market-objective_active_lim} 
\end{equation}
\begin{equation}
    \underline {g}_0^{Q} \leq {g}_{0t|ud}^{Q} \leq \overline {g}_0^{Q}, \label{eq:market-objective_reactive_lim}
\end{equation}
\begin{align} \text{variables } {}& \{ \varphi_{it|ud}, \xi_{it|ud}, f_{it|ud}^P, f_{it|ud}^Q, h_{0t|ud}^P, \nonumber \\ 
{}& h_{0t|ud}^Q, l_{it|ud}, v_{it|ud}, \hspace{10000pt minus 1fil} \nonumber \\ 
{}& i\in \mathcal{I}, n \in \Omega_i, u\in \mathcal{Y}^{pv}, d\in \mathcal{Y}^{de}, t\in \mathcal{T} \nonumber \},
\end{align}
\end{subequations}
where $\varpi >0$ and $\tau >0$ are additional prices for constraints violation; the associated slack variables $\xi_{it|ud}$ in (\ref{eq:market-objective_forward_limit}) and (\ref{eq:market-objective_backward_limit}) is about power flow limit, and $\varphi_{it|ud}$ in (\ref{eq:market-objective_vol_lim}) is voltage limit. These variables describe the amount of power that must be drawn from the system to avoid exceeding the maximum power flow limit and voltage limit. The variables are used to control the amount of energy drawn from the system and ensure that the system operates within safe limits.
\subsection{Coupling variables and algorithm process} \label{subsec:p2p_clearance_price}
Recall that utilize the bilateral P2P energy price $\lambda_{nmt|ud}$ in~(\ref{eq:sto_market_obj_fn}) as a coupling variable to connect Problems \textbf{P1} and \textbf{P2} as an iterative solution. This section proposes a bilateral P2P energy price $\lambda_{nmt|ud}^k$ at each iteration $k$ between \textbf{P1} and \textbf{P2}. Note that in describing \textbf{P1} and \textbf{P2}, we omitted the iteration index $k$ for notational simplicity, unless required. This price takes into account both the energy locational pricing (ELP), denoted by $\hat{\lambda}_{nmt|ud}$, and the main grid middle price (GMP), represented as $\lambda_{t}^{\text{mid}}$. By combining these two pricing models, we propose the integrated grid locational pricing (IGLP) step by step as follows.

In formulating the ELP, we introduce scenario-based pricing $\tilde{\lambda}_{nt|ud}^{k}$ for household $n \in \Omega_i$ at iteration $k$, expanding the formula in \cite{papavasiliou2017analysis} used for calculating distribution locational marginal pricing (DLMP) using the Lagrange multipliers of \textbf{P2},
\begin{equation}
\begin{aligned}
\tilde{\lambda}_{nt|ud}^{k}={}& C1_{t|ud}\gamma_{\mathcal{A}_{i}t|ud} + C2_{t|ud}\mu_{it|ud} + C3_{t|ud}\mu_{\mathcal{A}_{i}t|ud} \\
 {}& + C4_{t|ud}\eta_{it|ud}^{+} + C5_{t|ud}\eta_{it|ud}^{-}. \label{eq:main grid Cost Function}
\end{aligned}
\end{equation}
The parameters $C1_{t|ud}$, $C2_{t|ud}$, $C3_{t|ud}$, $C4_{t|ud}$, and $C5_{t|ud}$ can be determined as outlined in Appendix~\ref{appdix:DLMP_calulcation}. Moreover, considering that bilateral trade is reciprocal, ELP applies equally to both $n$ to $m$ trade and $m$ to $n$ trade. Therefore, ELP for a household $n$ trades with other households $m$ in distributed network at iteration $k$ can be calculated as
\begin{equation}
\label{eq:distributed_pricing}
\hat{\lambda}_{nmt|ud}^k = \hat{\lambda}_{nmt|ud}^{k-1} + \rho|\tilde{\lambda}_{nt|ud}^{k}-\tilde{\lambda}_{mt|ud}^{k}|, 
\end{equation}
where $\rho >0$ is a fine-tune parameter to reach the original solution. To formulate the GMP, we utilize the middle price between the buying price from the main grid $\overline {\lambda}_{t}$ and the selling price to the main grid $\underline {\lambda}_{t}$, expressed as
\begin{equation}
\lambda_{t}^{\text{mid}} = \frac{\overline \lambda_t + \underline \lambda_t}{2}. \label{eq:GMP_price}
\end{equation}
By utilizing the ELP in (\ref{eq:distributed_pricing}) and the GMP in (\ref{eq:GMP_price}), we propose the IGLP by
\begin{equation}
\label{eq:IGLP_price}
\lambda_{nmt|ud}^{k}
=\begin{cases}
\hat{\lambda}_{nmt|ud}^{k} - \lambda_{t}^{\text{mid}}, & \text{where $n$ sells to $m$} \\
\hat{\lambda}_{nmt|ud}^{k} + \lambda_{t}^{\text{mid}}, & \text{where $n$ buys from $m$}
\end{cases}
\end{equation}

Finally, the day-ahead market is summarized in \textbf{Algorithm~1}. Initially, the CO solves the economic dispatch problem \textbf{P1} at the community level and then transmits the net power of the community to the DSO. Subsequently, the DSO tackles the AC-OPF problem \textbf{P2} at the distribution level, updates coupling variables and sends back to the CO. The index $k$ captures the iteration between the CO and the DSO. The algorithm converges when the value of congestions remains unchanged in case congestion occurs. If there is no congestion, the algorithm converges to the solution with the smallest energy cost.
\begin{algorithm}[t]
\caption{Day-ahead electricity market}\label{alg:two-stage}
\begin{flushleft}
        \parbox[t]{\dimexpr\linewidth-\algorithmicindent\relax}{%
        \setlength{\hangindent}{0pt}%
        \textbf{Input:} Predicted PV generation and demand scenario of households;\\
        \textbf{Output:} ESS strategies $P_{nt}^{bc}, P_{nt}^{bd}$;\\
        \textbf{Initialization:} $k=1$, $\hat{\lambda}_{nmt|ud}^k=0$, $\tilde{\lambda}_{mnt|ud}^k=0$, and given $\overline \lambda_{t}$ and $\underline \lambda_{t}$ from main grid;}\strut
\end{flushleft}
\begin{algorithmic}[1]
\While{true}
\State \parbox[t]{\dimexpr\linewidth-\algorithmicindent\relax}{%
    \setlength{\hangindent}{0pt}%
      CO solves the stochastic optimization in Problem \textbf{P1} at community level;}\strut
\State \parbox[t]{\dimexpr\linewidth-\algorithmicindent\relax}{%
    \setlength{\hangindent}{0pt}%
      CO sends $s_{nmt|ud}, b_{nmt|ud}, s_{n0t|ud},$ and $b_{0nt|ud}$ to DSO;}\strut
\State \parbox[t]{\dimexpr\linewidth-\algorithmicindent\relax}{%
    \setlength{\hangindent}{0pt}%
      DSO solves cost minimization Problem in \textbf{P2} at network distribution level;}\strut
\State \parbox[t]{\dimexpr\linewidth-\algorithmicindent\relax}{%
    \setlength{\hangindent}{0pt}%
      DSO updates IGLP $\lambda_{nmt|ud} = \lambda_{nmt|ud}^k$ according to (\ref{eq:IGLP_price}) then sends back to CO;}\strut
\If{stop condition is met} 
    \State CO sends the ESS strategies $P_{nt}^{bc}$, 
    \State $ P_{nt}^{bd}$ to household $n$;
    \State return;
\EndIf
\State $k \leftarrow k+1$;
\EndWhile
\end{algorithmic}
\end{algorithm}
\section{Energy market clearing problem}\label{sec:Energy-market-clearing}
Every hour in one day presents new decision-making challenges for both the CO and DSO, as depicted in Fig.~\ref{fig:overview}. The proposed method addresses these challenges by applying day-ahead solutions to real-time operations when true data on PV generation and demand becomes available, i.e., we only need to consider the realized scenarios $u\in \mathcal{Y}^{pv}$ and $d\in \mathcal{Y}^{de}$, and  thus, for notational simplicity we exclude the subscripts~$_{|ud}$ hereafter. Consequently, this section explores how the proposed approach resolves issues in the real-time market and subsequently calculates energy market payments based on power exchange within the power grid.
\subsection{Real-time operation} \label{subsec:infeasible}
Once real-time data on PV generation $G_{nt}^{pv}$ and demand $L_{nt}^{de}$ are available, the energy market faces a new problem to determine the optimal strategies of households, considering ESS, bilateral trades with the grid, and transactions between households under real power exchange in the distribution network. This decision-making process entails solving a deterministic problem based on PV generation and demand, then energy balance constraints (\ref{eq:house_net_active_power})--(\ref{eq:house_net_reactive_power}) are revised as
\begin{subequations}
\begin{equation}
    P_{nt}^{P} = G_{nt}^{pv} + P_{nt}^{bd} + e_{nt}^{bc} - L_{nt}^{de} - P_{nt}^{bc} - e_{nt}^{bd}, \label{eq:RT_house_net_active_power}
\end{equation}
\begin{equation}
    P_{nt}^{P} = s_{n0t} + \sum_{m\in \omega_n} s_{nmt} - b_{n0t} - \sum_{m\in \omega_n} b_{nmt}, \label{eq:RT_house_energy_balance}
\end{equation}
\begin{equation}
    P_{nt}^{Q} = - L_{nt}^{Q}, \label{eq:RT_house_net_reactive_power}
\end{equation}
\end{subequations}
where $e_{nt}^{bc}$ and $e_{nt}^{bd}$ facilitate the adjustment of ESSs patterns in real-time operation to prevent constraint violations induced by forecasting errors, then the model of battery (\ref{eq:battery_SOC})--(\ref{eq:battery_binary}) is updated as
\begin{subequations}
\label{eq:RT_battery}
\begin{align}
E_{nt+1}^{RT}= E_{nt}^{RT}+\biggr[\eta (P_{nt}^{bc}+e_{nt}^{bc}) - (P_{nt}^{bd} + e_{nt}^{bd})/\eta \biggr]\bigtriangleup_t, \label{eq:RT_battery_SOC}
\end{align}
\begin{equation}
\underline {E}_n \leq E_{nt+1}^{RT} \leq \overline {E}_n, \label{eq:RT_battery_SOC_limit}
\end{equation}
\begin{equation}
0 \leq P_{nt}^{bc} + e_{nt}^{bc} \leq \psi_{nt} P^b_n, \label{eq:RT_bc_limit}
\end{equation}
\begin{equation}
0 \leq P_{nt}^{bd} + e_{nt}^{bd} \leq (1-\psi_{nt}) P^b_n, \label{eq:RT_bd_limit}
\end{equation}
\begin{equation}
\psi_{nt} \in [0,1],  \label{eq:RT_battery_binary}
\end{equation}
\end{subequations}
At time slot $t+1$, the amount of energy of household $n$, the communities have to be adjusted, denoted by $\bigtriangleup E_{nt+1}$, can be calculated by the stored energy for day-ahead $E_{nt+1}^{DA}$ and the stored energy for real-time operation $E_{nt+1}^{RT}$ such as
\begin{equation}
    \bigtriangleup E_{nt+1} = E_{nt+1}^{DA} - E_{nt+1}^{RT}, \label{eq:deviation_SOC}
\end{equation}
where $\bigtriangleup E_{nt+1}$ takes positive value when the CO charges energy to battery of a specific household $n$ to alleviate congestion; otherwise, negative value means the community discharges energy from the battery of household $n$.

In this regard, we reformulate the CO and DSO interaction at time slot $t$ into a single optimization problem \textbf{P3} as follows:

\vspace{2mm}
\noindent \textbf{P3: Real-time clearing problem
}
\begin{align}
\text{minimize ~~} {}& \sum_{n\in\mathcal{N}} (\overline \lambda_{t} b_{n0t} - \underline \lambda_{t} s_{n0t} + \overline \lambda_{t}{g}_{0t}^{P} + \theta|\bigtriangleup E_{nt+1}| ) \label{eq:RT-SWM}\\
\text{subject to ~} {}& \text{ESS: (\ref{eq:RT_battery_SOC})--(\ref{eq:RT_battery_binary}), (\ref{eq:deviation_SOC})}, \hspace{10000pt minus 1fil} \nonumber \\
&\text{P2P:  (\ref{eq:net_trade_p2p})--(\ref{eq:bilateral_p2p}), (\ref{eq:RT_house_net_active_power})--(\ref{eq:RT_house_net_reactive_power})}, \nonumber \\
{}& \text{NET: (\ref{eq:ac-opf_forward_limit})--(\ref{eq:ac-opf_vol_lim}), (\ref{eq:ac-opf_net_active_flow})--(\ref{eq:ac-opf_net_reactive_flow})}, \nonumber \\
\text{variables ~~} {}&  \{b_{nmt}, s_{nmt}, T_{nmt},P^P_{nt}, b_{n0t}, s_{n0t},\nonumber \\ 
{}& v_{it}, l_{it}, h_{it}^P, h_{it}^Q, e_{nt}^{bc}, e_{nt}^{bd},  \nonumber \\
{}& n\in \mathcal{N}, i\in \mathcal{I},n\in\Omega_i, m\in\omega_n, t\in \mathcal{T} \nonumber \},
\end{align}
where $\theta >0$ is a penalty parameter to adjust the solution in real-time operation at time slot $t$. Note that problem \textbf{P3} solves the AC-OPF problem of the DSO. Therefore, the bilateral P2P energy price $\lambda_{nmt}$ is determined after solving \textbf{P3} and will be utilized in the next subsection.
\begin{figure}
\centering
\includegraphics[width=\columnwidth]{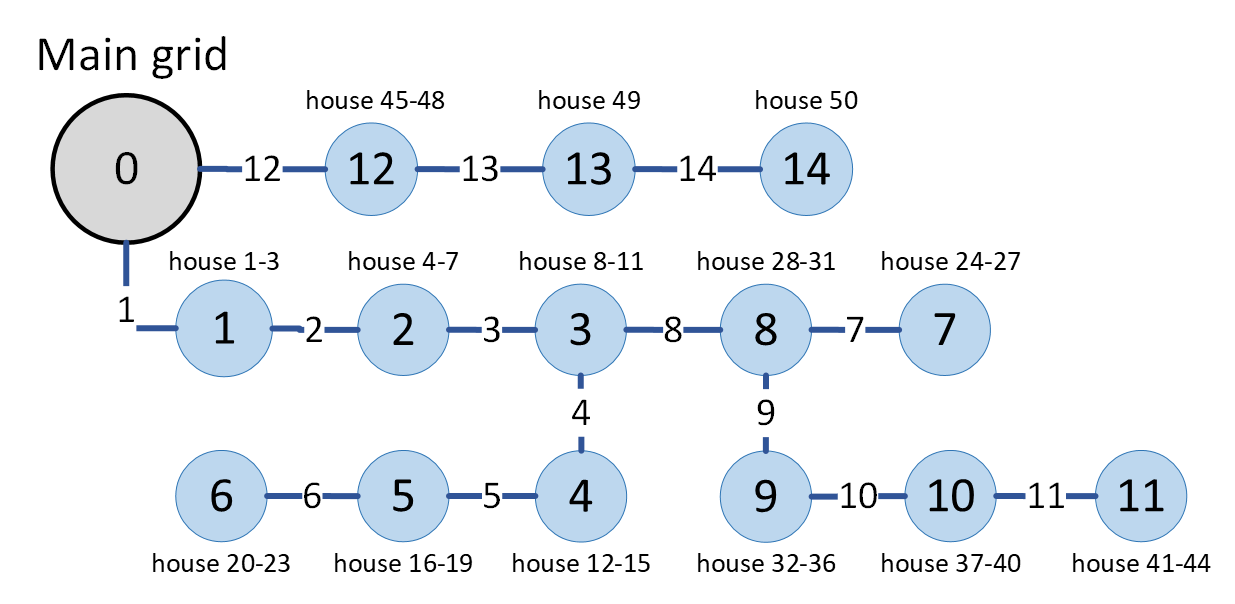}
\caption{Distributed network test system}
\label{fig:network_diagram}
\end{figure}
\begin{figure}[t]
\centering
\includegraphics[width=\columnwidth]{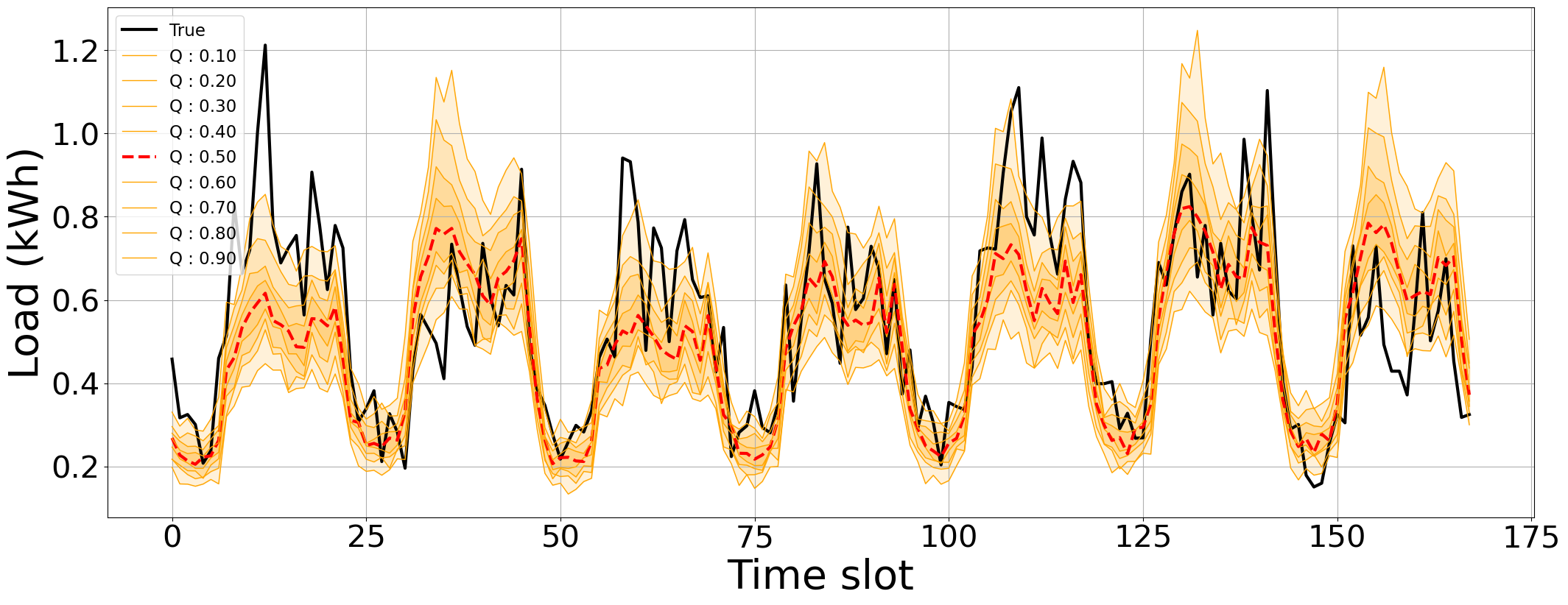}
\includegraphics[width=0.98\columnwidth]{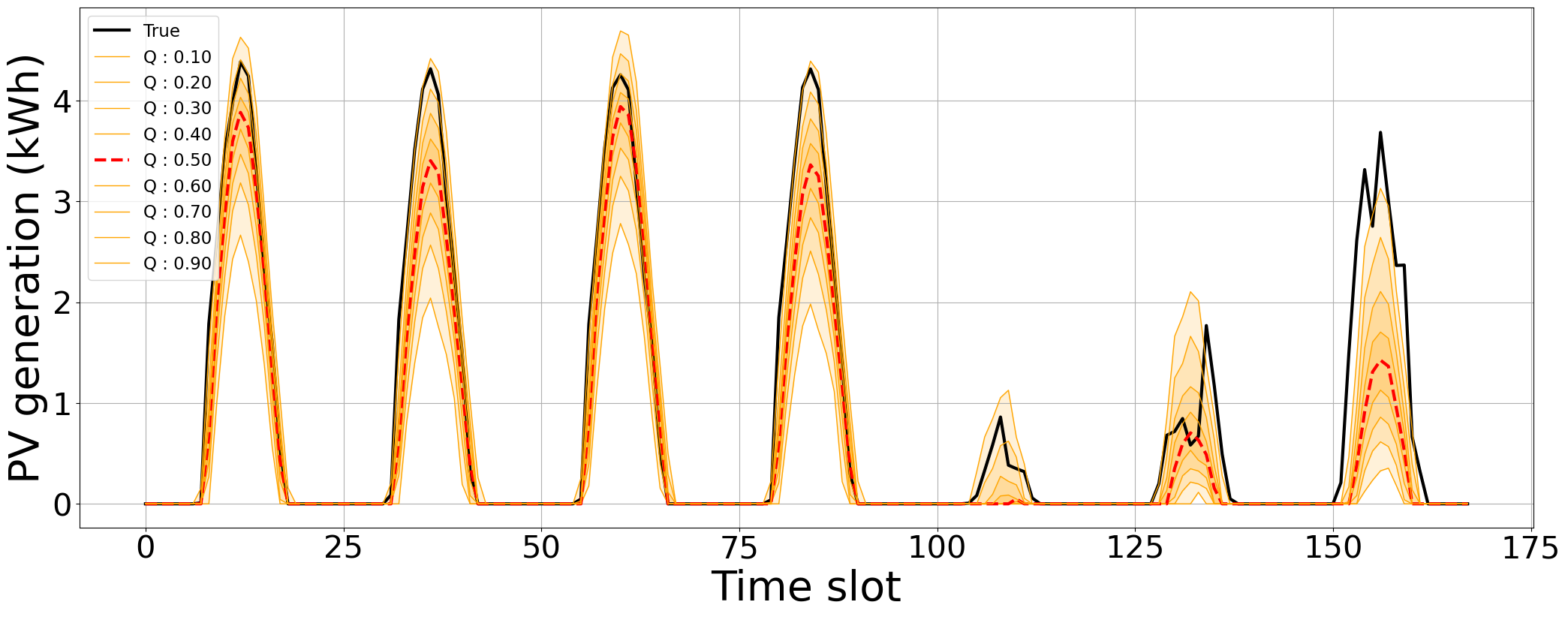}
\caption{Example of quantile-based demand/PV scenario for a household. The quantile forecasting algorithm generates the median value (red) and each quantile value (orange). The true demand/PV energy is illustrated by the black line.}
\label{fig:PV_demand_quantile}
\end{figure}
\subsection{Payment} \label{subsec:settlement}
While the preceding section addresses the issues within the energy market, it does not accurately reflect the true energy costs incurred by households. Hence, this subsection will delve into the actual energy consumption costs and the selling profit associated with the energy that a household $n$ pays or receives at time slot $t$. The true energy costs for both the community and households can be computed at the end of the day throughout each hour based on the recorded transitive energy in smart meters. Consequently, the cost of the community at time slot $t$, denoted by $C_{\mathcal{N}t}$, can be reformulated as follows:
\begin{align}
    C_{\mathcal{N}t} =& \sum_{n\in \mathcal{N}} C_{nt}, \\
    C_{nt} =& \overline {\lambda}_{t}b_{n0t} - \underline {\lambda}_{t}s_{n0t} - \sum_{m\in \omega_n}\lambda_{nmt}T_{nmt}.
\end{align}
The cost for an individual household $n$ is represented as $C_{nt}$, and the value negative indicates that household $n$ earns money by selling its energy. The actual energy exchanged with the main grid at time slot $t$ is denoted as $b_{n0t}$ and $s_{n0t}$, which are recorded in the smart meter for household $n$. The IGLP, $\lambda_{nmt}$, which is determined by solving Problem \textbf{P3}, and $T_{nmt}$ reflects the actual net power of P2P trades when household $n$ buys/sells from/to other household $m$. Note that this true energy costs will be used to show efficiency of the proposed method in the next section.
\section{Performance Evaluation}\label{sec:performance}
We evaluate the performance of the proposed method and compare it with both an \textit{offline} approach and a \textit{point-forecasting} method \cite{van2020integrated}. The offline method achieves the minimum cost for day-ahead energy scheduling because it assumes the perfect knowledge of future. By contrast, with point-forecasting method in \cite{van2020integrated}, energy market is solved based on the point forecasting data. The results of these methods will be compared with those of the proposed method, all based on the Gurobi optimizer \cite{gurobi2021gurobi}.
\begin{figure}
\centering
\includegraphics[width=\columnwidth]{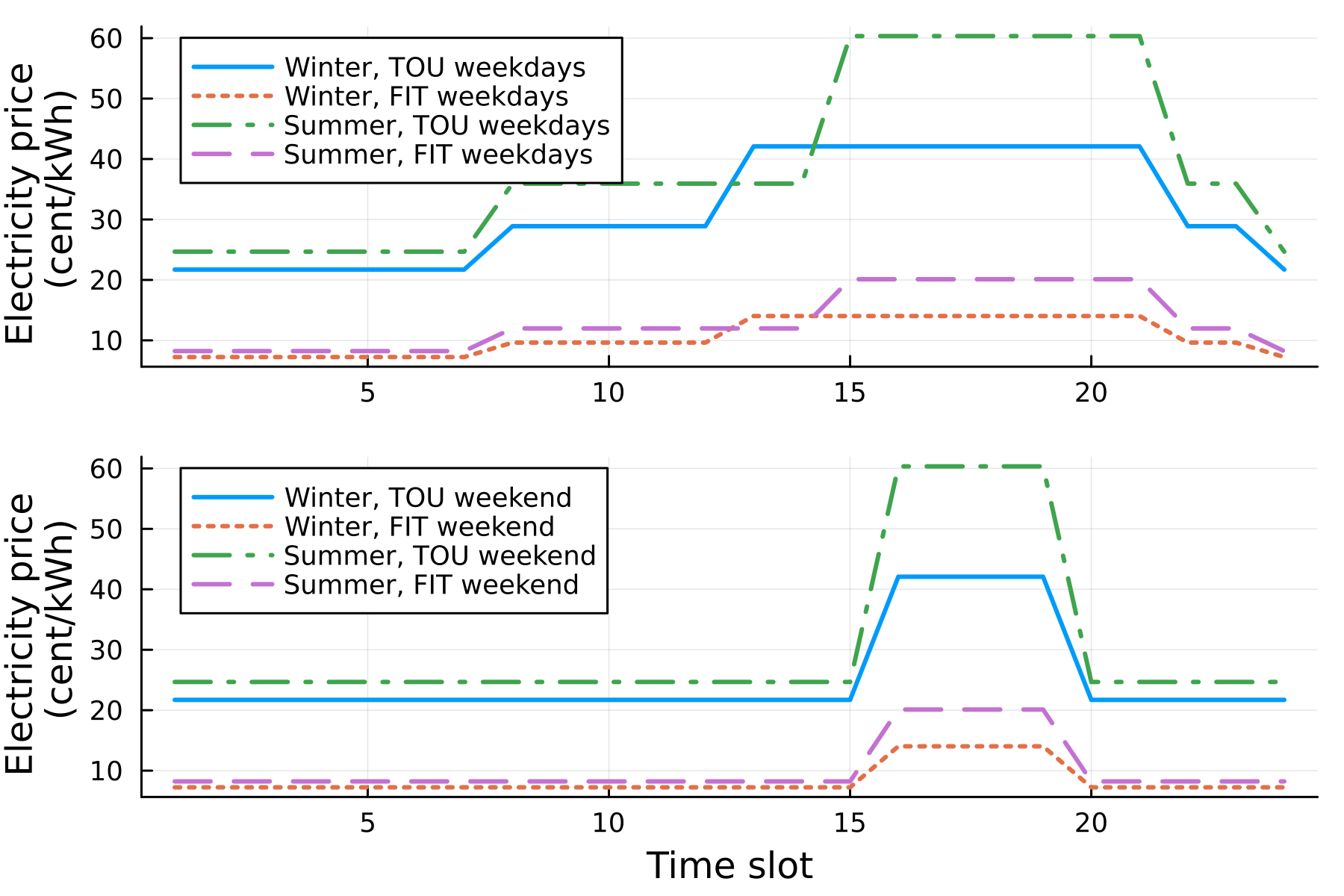}
\caption{Electricity price of main grid during weekdays (top) and weekend (bottom).}
\label{fig:electrical_price}
\end{figure}
\begin{table}
\vspace{0.5em}
\caption{Battery setup for a household $n$.}
\label{tab:battery_setup}
\centering
\resizebox{0.5\columnwidth}{!}{%
\begin{tabular}{|c|c|}
\hline
Parameters          & Values \\ \hline
$\eta$ &    95(\%)    \\ \hline
$\overline {E}_n$                 &    13.5 (kWh)   \\ \hline
$\underline {E}_n$              &     0.64 (kWh)  \\ \hline
$P^b_n$                  &    3.3 (kW)   \\ \hline
$ {E}_{n1} = {E}_{n24}$                  &    0.64 (kWh)   \\ \hline
\end{tabular}}
\end{table}

A case study is based on the modified IEEE 15-bus distribution system \cite{kim2019p2p}. The network comprises two years of 50 real household datasets, provided by LG Electronics after anonymization, featuring individual PVs, ESS, household demands, and inter-household trading along with their partners and the upstream main grid. The geographical distribution of these households is visualized in Fig.~\ref{fig:network_diagram}, while the specific configuration of the energy storage system is illustrated in Table~\ref{tab:battery_setup}. The parameters for the energy market mechanism is in Table~\ref{tab:addition_parameters}. An example of electrical prices during one week is depicted in Fig.~\ref{fig:electrical_price}. These prices, relevant for buying energy from grid time-of-use (TOU) and selling energy to grid feed-in tariffs (FIT), dynamically change based on the time of day, the day of the week, and whether it is a weekday or weekend. We assume 10\% of reactive power to active power and the time slot duration $\bigtriangleup_t = 1$. In forecasting hourly energy generation and consumption within a household in the community is illustrated in Fig.~\ref{fig:PV_demand_quantile}. By employing TransLSTM model \cite{song2022dttrans}, we generate diverse scenarios for both PV generation and demand encompassing quantiles ranging from $0.1$ to $0.9$.
\subsection{Profit evaluation and comparison}
\begin{figure}[t]
     \centering
     \begin{subfigure}[b]{0.48\textwidth}
         \centering
         \includegraphics[width=\textwidth]{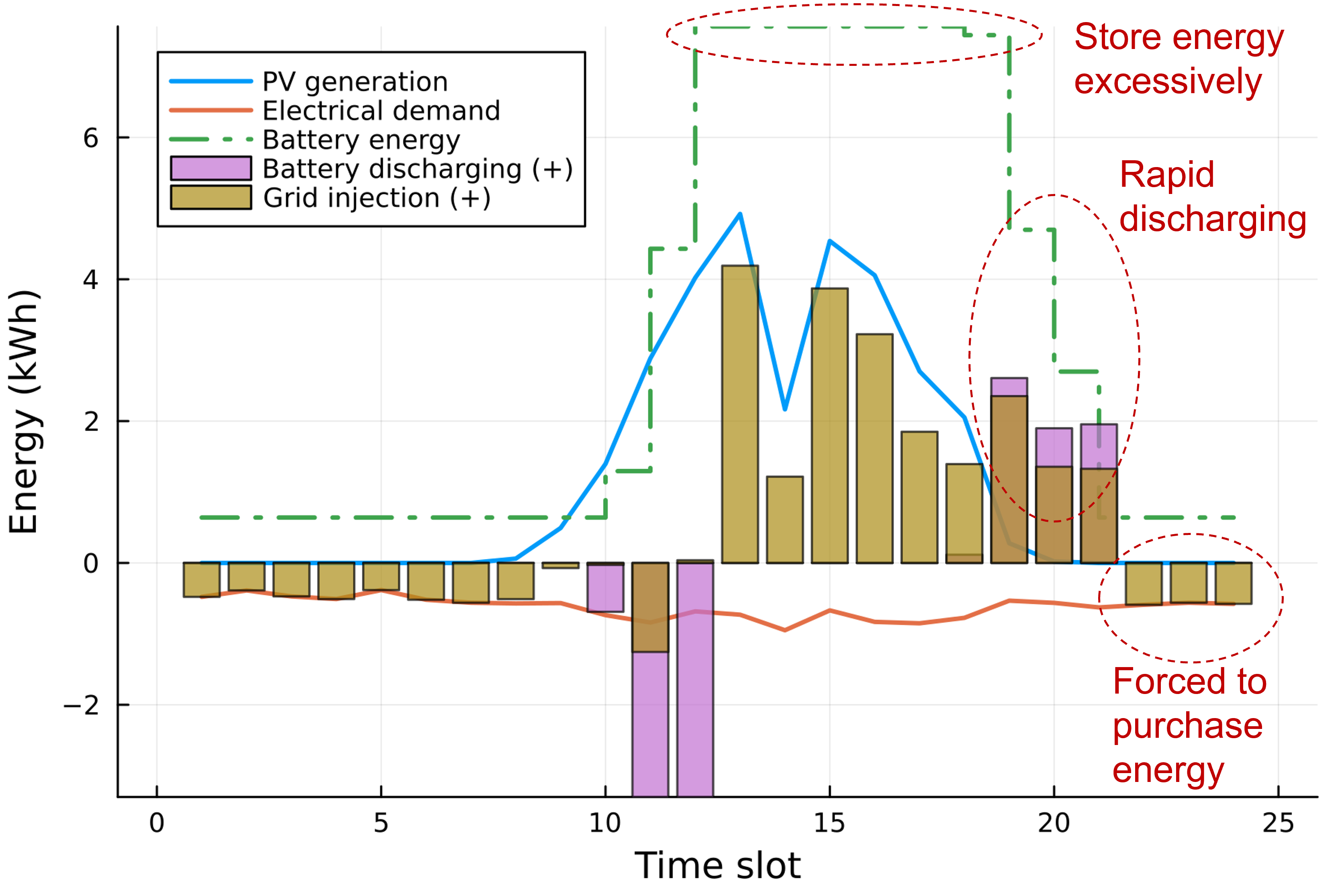}
         \caption{Point-forecasting method in \cite{van2020integrated}}
         \label{fig:cost_deter_method}
     \end{subfigure}
     \hfill
     \begin{subfigure}[b]{0.48\textwidth}
         \centering
         \includegraphics[width=\textwidth]{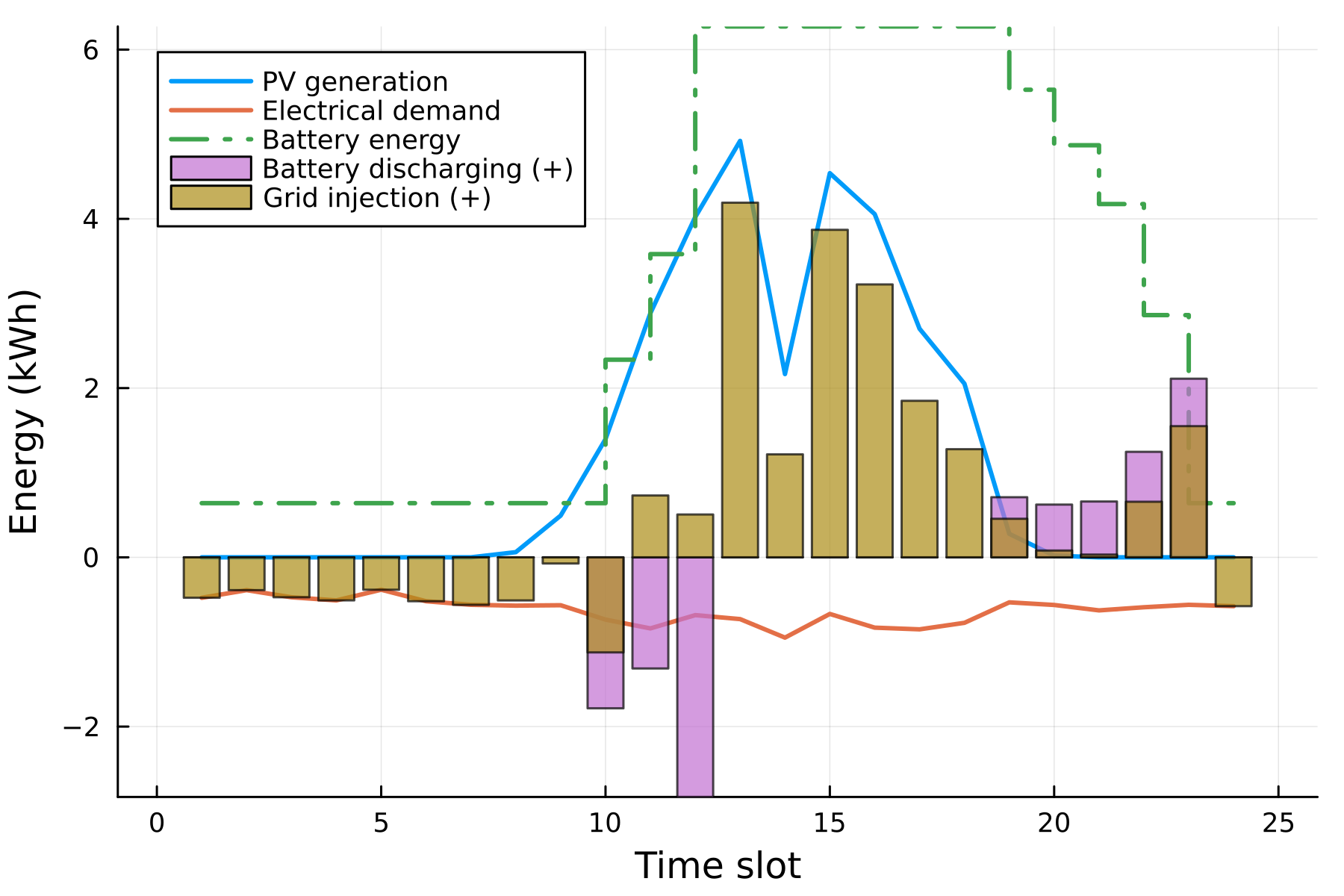}
         \caption{Proposed method}
         \label{fig:cost_proposed_method}
     \end{subfigure}
        \caption{Energy scheduling of household indexed 10 (positive values are for PV generation, injection into the grid and battery discharging).}
        \label{fig:energy_scheduling_aday}
\end{figure}
To visualize energy scheduling of a typical household, we randomly select a household and present the results for one day, in Fig.~\ref{fig:energy_scheduling_aday}. For this particular household, both the point-forecasting and the proposed method purchase energy to meet immediate demand and refrain from acquiring energy for battery charging during time slots 1 to 10. Between time slots 10 and 12, both methods attempt to charge the PV generation into the ESS. However, the stored energy differ between the two methods; the point-forecasting method stores energy in the ESS excessively. In next time slot from time slots 13 to 18, the ESS strategies remain the same in two methods, in which all surplus PV generation is utilized for selling to the main grid and other households, as well as to meet their own demand. However, the ESS strategies of two methods differ in the remaining time slots from 19 to 24. In the point-forecasting method, the discharge is rapid from time slot 19 to 21, and in return energy is purchased from the grid during the last three time slots, as shown in Fig.~\ref{fig:cost_deter_method}. By constrast, the proposed method discharges slowly and only makes a purchase during the last time slot. In this example, it can be seen that the proposed ESS strategies conserve energy for future use and promote the sale of surplus energy when there is no sunlight, as illustrated in Fig.~\ref{fig:cost_proposed_method}. Thanks to this, the proposed method reduces dependence on purchasing energy from the main grid and increases profit by selling energy more efficiently than the point-forecasting method.
\begin{table}
\vspace{0.5em}
\caption{Parameter setup for energy market mechanism.}
\label{tab:addition_parameters}
\centering
\resizebox{0.5\columnwidth}{!}{%
\begin{tabular}{|c|c|}
\hline
Parameters          & Values \\ \hline
$\varpi$ in Problem \textbf{P1} &    10    \\ \hline
$\tau$ in Problem \textbf{P2}                &    10   \\ \hline
$\rho$                &    1   \\ \hline
$\theta$ in Problem \textbf{P3}             &   50  \\ \hline
\end{tabular}}
\end{table}
\begin{table}
\vspace{0.5em}
\caption{Profit evaluation of three different methods.}
\label{tab:energy_costs_three_methods}
\centering
\resizebox{\columnwidth}{!}{%
\begin{tabular}{|l|ll|}
\hline
\multirow{2}{*}{}        & \multicolumn{2}{c|}{ \textbf{ \makecell{Total community profit\\in a week}}}                              \\ \cline{2-3} 
                         & \multicolumn{1}{c|}{Summer season} & \multicolumn{1}{c|}{Winter season} \\ \hline
Point-forecasting method \cite{van2020integrated} & \multicolumn{1}{l|}{80810.9}      &  50832.2                          \\ \hline
Proposed method          & \multicolumn{1}{l|}{83994.2 }      & 55405.8                           \\ \hline
Offline method           & \multicolumn{1}{l|}{86155.7 }      & 57661.9                           \\ \hline
\end{tabular}}
\end{table}

As illustrated in Table~\ref{tab:energy_costs_three_methods}, it becomes evident that the proposed method is more efficient than the point-forecasting approach during a week. It should be noted that Problems \textbf{P1}, \textbf{P2}, and \textbf{P3} are geared toward minimizing community energy costs. By selling energy, costs become negative, and thus become profits as depicted in Table~\ref{tab:energy_costs_three_methods}. Interestingly, the proposed method reduces the optimality gap by 60\% and 67\% in summer and winter, respectively. The optimality gap in this paper refers to the additional cost compared to that of the offline method having the perfect knowledge of future. In their efforts to maximize energy profits, households aim to maximize the self-consumption of PV generation and engage in energy transactions with other households. Alternatively, they may seek to sell surplus energy to maximize profits instead of relying solely on the main grid. Furthermore, our method facilitates informed decision-making amid generation and demand uncertainty and variability. Consequently, the energy profits in our method are higher than those of the point-forecasting method.
\begin{figure}
\vspace{0.5em}
\centering
\includegraphics[width=\columnwidth]{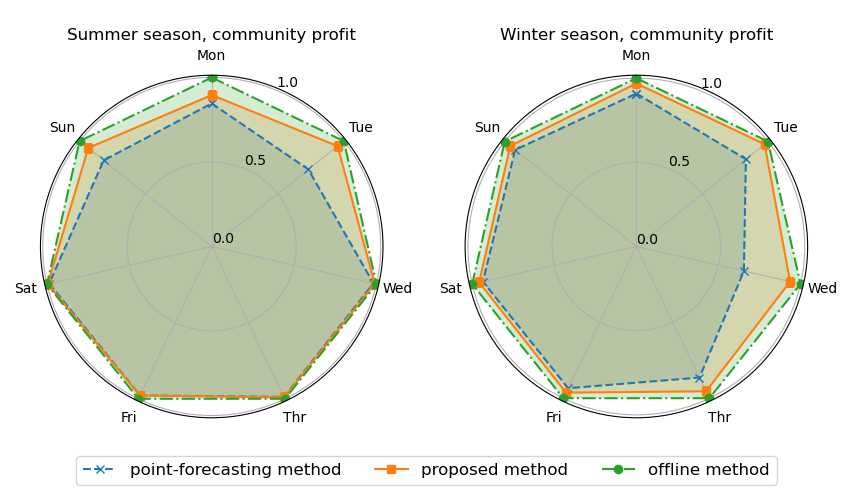}
\caption{A comparison of profit of three methods during a week. Each profit is normalized by the \textit{offline} method's profit.}
\label{fig:daily_energy_costs_three_methods}
\end{figure}

Furthermore, energy profits vary across different seasons of the year. In summer, community profit is higher than in the winter season, meaning there is 30\% higher energy costs during the winter season compared to the summer season, as depicted in Table~\ref{tab:energy_costs_three_methods}. The higher energy costs in winter can be attributed to several factors. Firstly, the reduced sun generation energy leads to a decrease in the amount of solar power that can be harvested, increasing dependence on the main grid for energy and consequently raising costs. Moreover, the higher prices of electricity in the winter season further contribute to the higher energy costs, as shown in Fig.~\ref{fig:electrical_price}. Additionally, cold weather and increased heating demands from households further strain the energy supply, leading to higher energy expenses.
\begin{table}
\vspace{0.5em}
\caption{Analyzing the constrains of violation for 24 hours.}
\label{tab:analyzing_constrains_week}
\centering
\resizebox{\columnwidth}{!}{%
\begin{tabular}{|l|l|}
\hline
                               & \textbf{Constraint violations}                                                                                                                                       \\ \hline
\makecell{ Point-forecasting \\ method \cite{van2020integrated}} & \begin{tabular}[c]{@{}l@{}}
Line 2 from bus 1 to 2 overflow at time slot 14\\
Line 2 from bus 1 to 2 overflow at time slot 15\\
Line 2 from bus 1 to 2 overflow at time slot 16\\
Line 2 from bus 1 to 2 overflow at time slot 17\\\end{tabular} \\ \hline
Proposed method               & Constraints are all met at all time slots                                                                                                                                                            \\ \hline
\end{tabular}}
\end{table}

Moreover, for a comprehensive evaluation of performance over a week, as shown in Table~\ref{tab:energy_costs_three_methods}, we present a detailed breakdown of community energy profits for a single day in Fig.~\ref{fig:daily_energy_costs_three_methods}. The results consistently demonstrate the superior performance of the proposed method over the point-forecasting method across each day of the assessment. This underscores a higher level of reliability and accuracy in estimating community energy profits with the proposed method compared to the point-forecasting method.
\begin{figure}[t]
\centering
\begin{subfigure}[b]{0.48\textwidth}
 \centering
 \includegraphics[width=\textwidth]{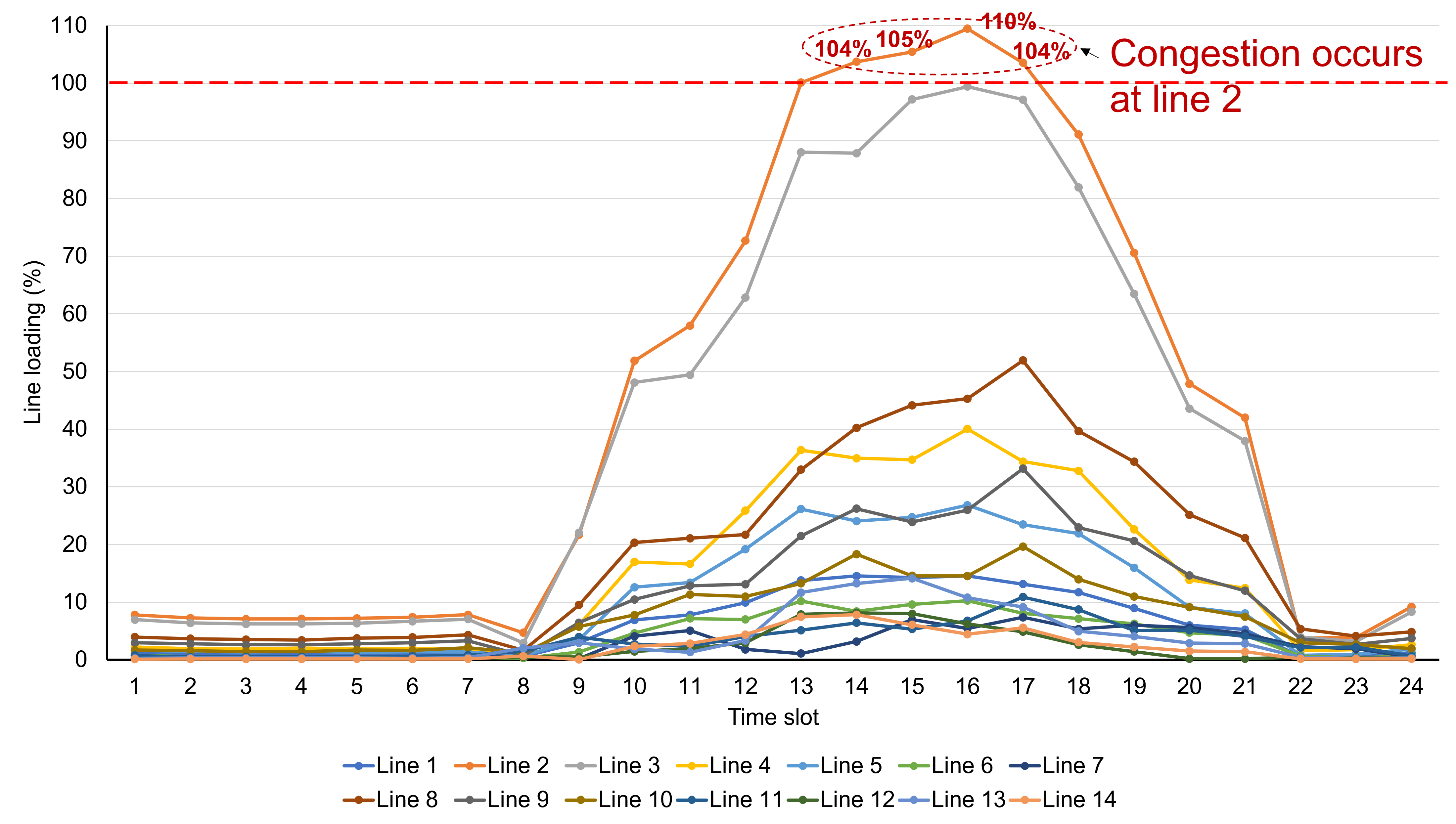}
 \caption{Point-forecasting method in \cite{van2020integrated}}
 \label{fig:pf-point-forecasting}
\end{subfigure}
\hfill
\begin{subfigure}[b]{0.48\textwidth}
 \centering
 \includegraphics[width=\textwidth]{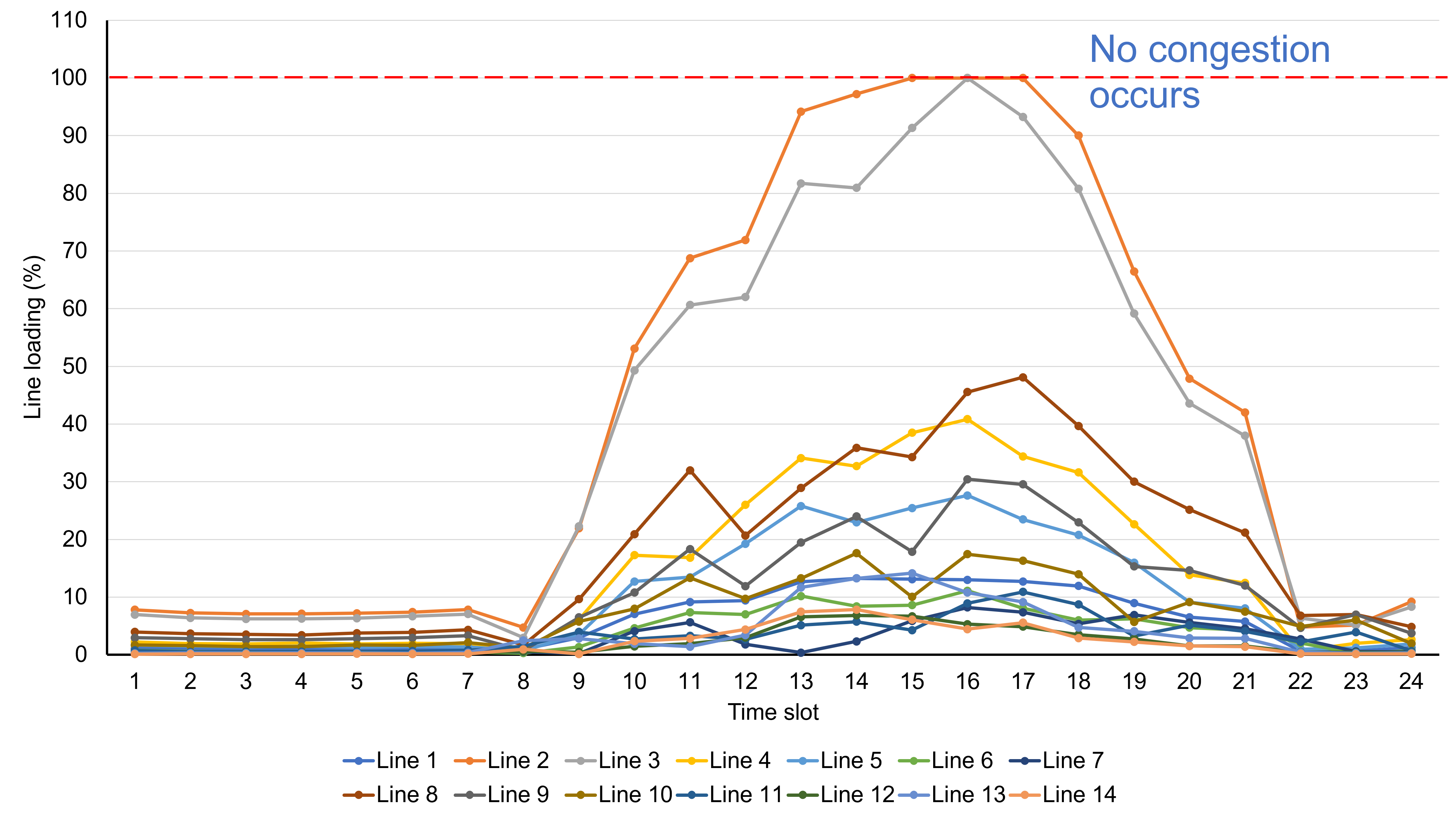}
 \caption{Proposed method}
 \label{fig:pf-proposed-method}
\end{subfigure}
\caption{Power flows at all buses over a 24-hour period.}
\label{fig:Power_flows_all_buses}
\end{figure}
\subsection{Congestion mitigation} 
To demonstrate the impact of uncertain PV generation-demand on the maximum power flow constraints, all constraints are evaluated across all time slots, and the outcomes of interest summarized are detailed in Table~\ref{tab:analyzing_constrains_week}. The proposed method consistently satisfies all the constraints for every line throughout all time slots. By contrast, the point-forecasting method struggles to meet the constraints of power flow and voltage and exceeds the congestion constraint in the afternoon when PV generation is high.

In the proposed method, the key lies in addressing congestion proactively as a forward-looking challenge; ESS strategies focus on balancing PV generation and demand as well as managing the power flow in the distribution network over a 24-hour period, as shown in Fig.~\ref{fig:Power_flows_all_buses}. The point-forecasting method engages in excessive power exchange, notably from 14:00 to 18:00, leading to congestion occurrences depicted in Fig.~\ref{fig:pf-point-forecasting}. By contrast, the proposed ESS strategies follow a distinctive pattern to mitigate uncertainties in PV generation and demand. They aim to decrease peak net injection, flatten energy net injection curves, prioritize enhancing self-consumption, and reduce power system exchanges during time slots marked by heightened congestion rates, as illustrated in Fig.~\ref{fig:pf-proposed-method}.

This comparative analysis underscores the effectiveness of the proposed approach in transitioning net injection from peak to average power exchange. Additionally, it emphasizes the proposed method's capacity to provide a more resilient solution capable of responding effectively to changes in the power system, thereby mitigating congestion.

\begin{figure}[t]
\centering
\includegraphics[width=\columnwidth]{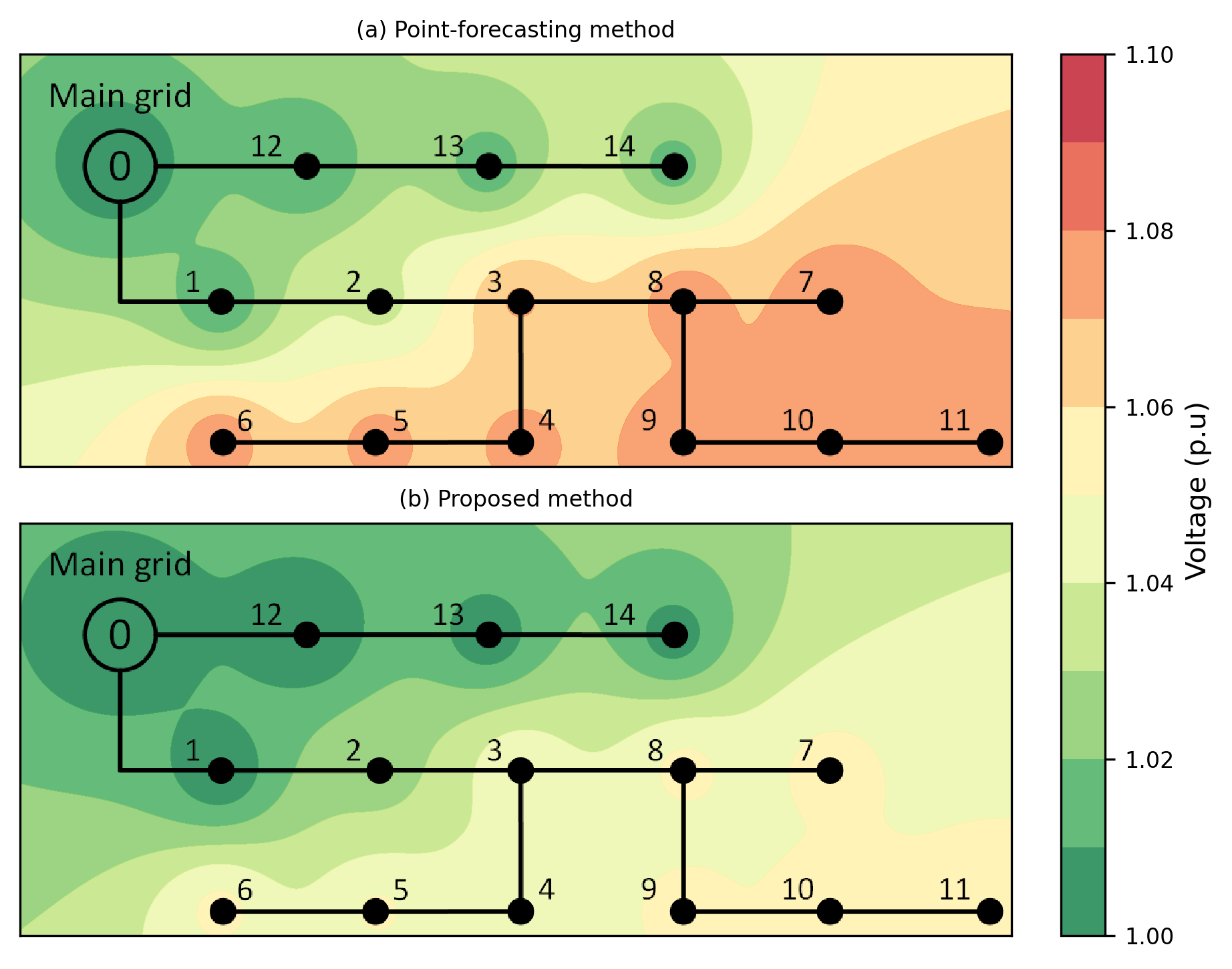}
\caption{The voltage on all buses during time slot 14 when the network is congested.}
\label{fig:voltage-slot-14}
\end{figure}

\begin{figure}[t]
\centering
\includegraphics[width=\columnwidth]{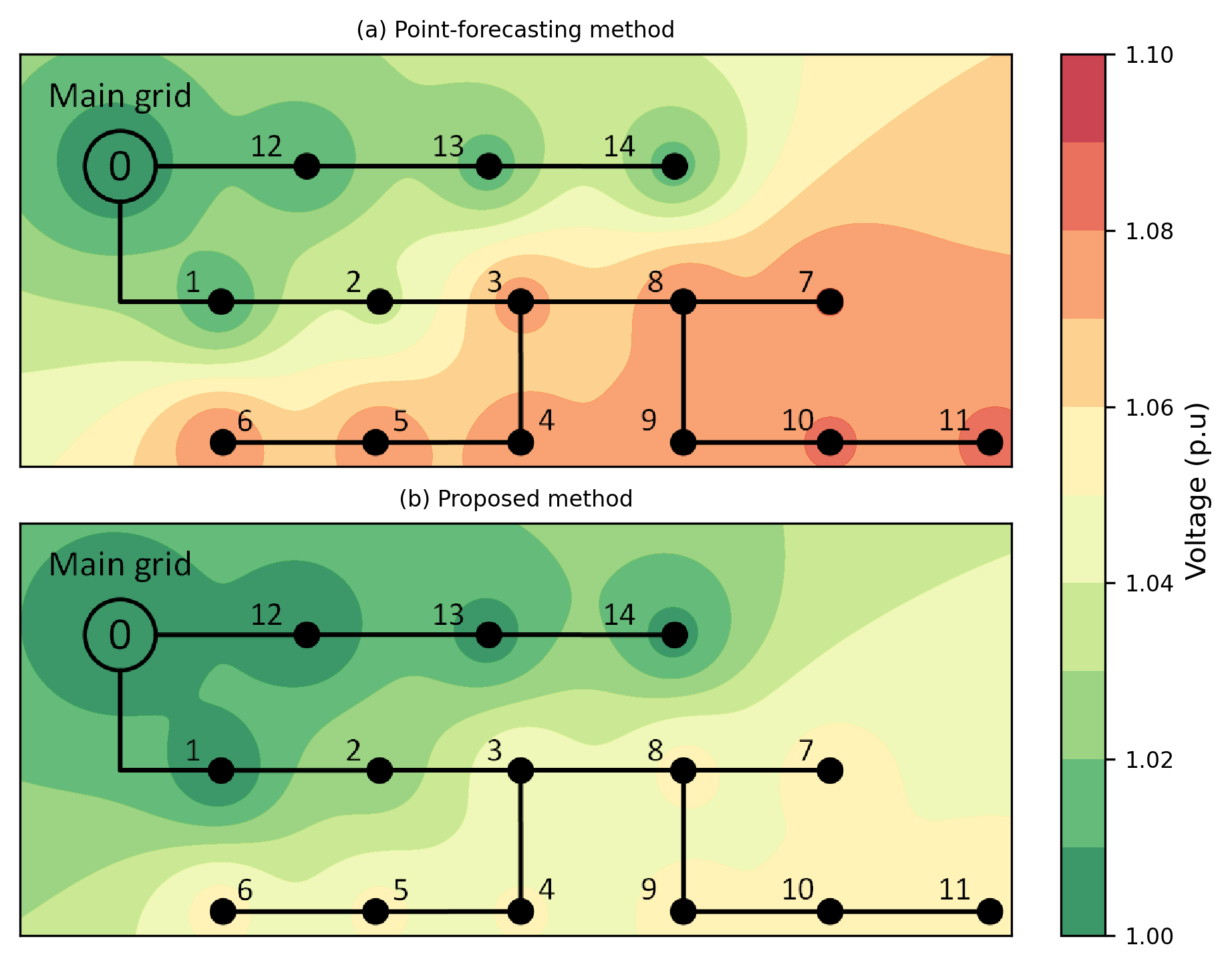}
\caption{The voltage on all buses during time slot 15 when the network is congested.}
\label{fig:voltage-slot-15}
\end{figure}

\begin{figure}[t]
\centering
\includegraphics[width=\columnwidth]{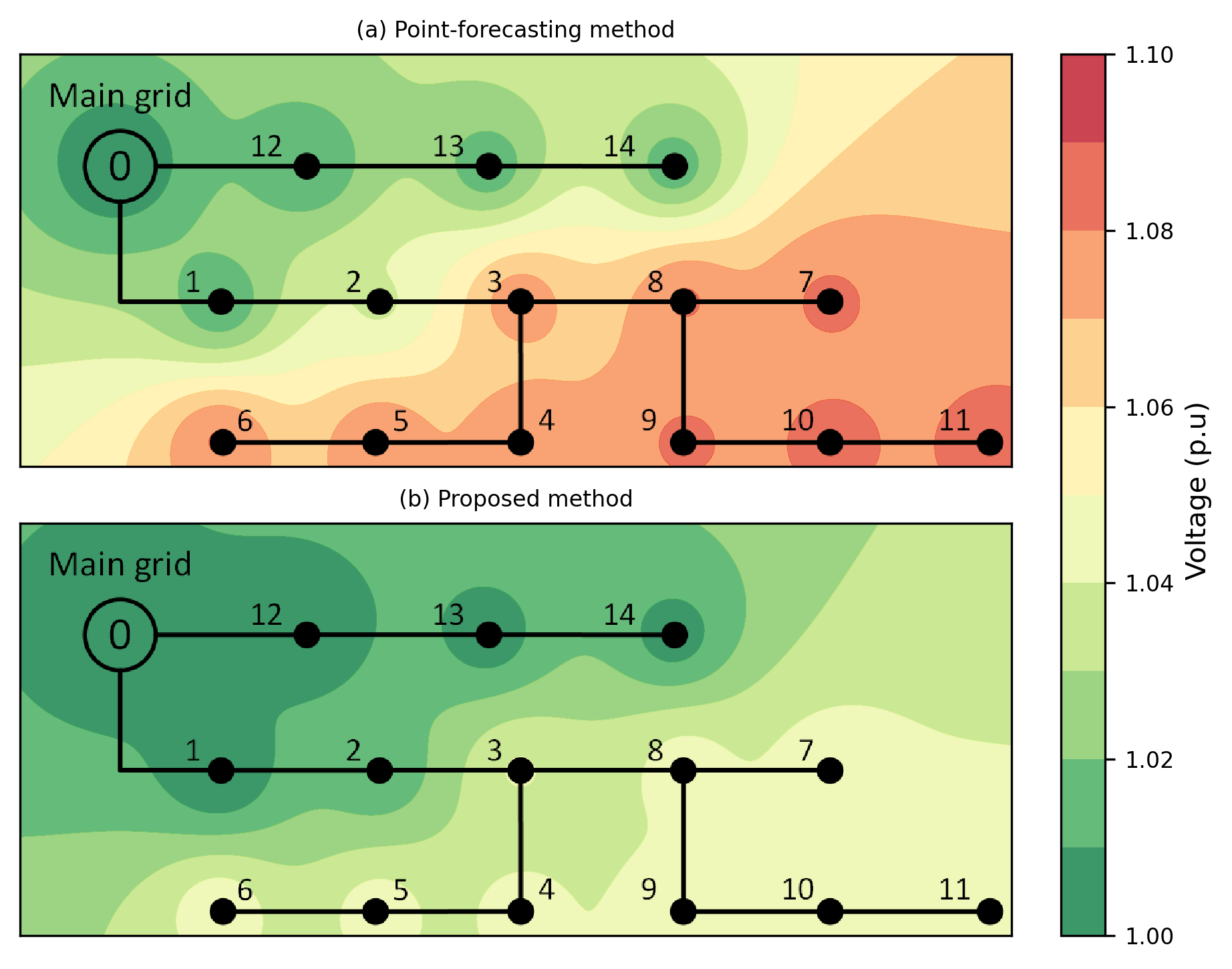}
\caption{The voltage on all buses during time slot 16 when the network is congested.}
\label{fig:voltage-slot-16}
\end{figure}

\begin{figure}
\centering
\includegraphics[width=\columnwidth]{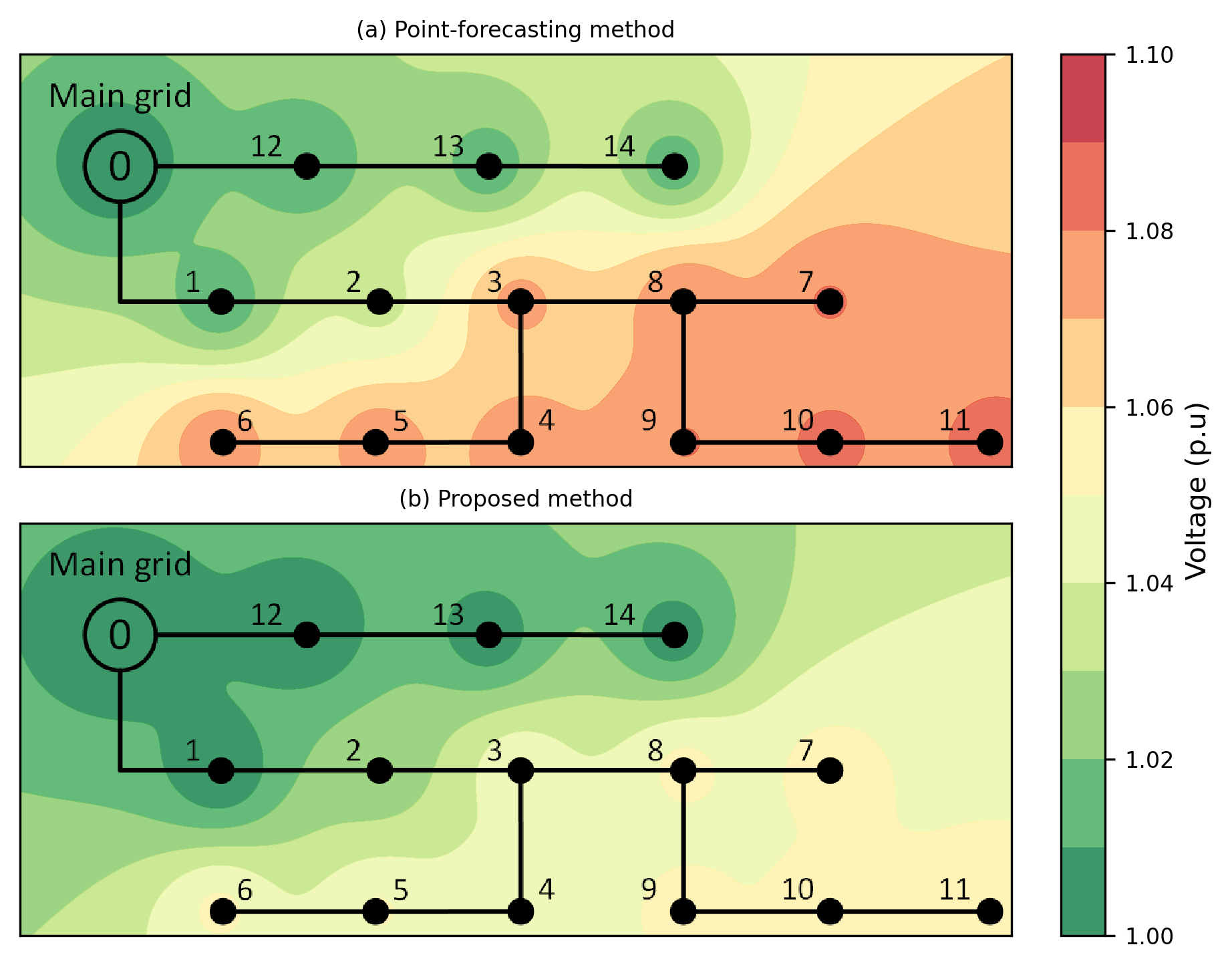}
\caption{The voltage on all buses during time slot 17 when the network is congested.}
\label{fig:voltage-slot-17}
\end{figure}

Finally, we investigate voltage violation and its mitigation. Figs.~\ref{fig:voltage-slot-14} to~\ref{fig:voltage-slot-17} analyze the voltage across all buses in the congested time slots of Fig.~\ref{fig:Power_flows_all_buses}. The results indicate that the proposed method successfully mitigates voltage violation during the congested time slots from 14 to 17 and increase the voltage stability margin of the network. In summary, the proposed method offers various advantages, including an elevated level of accuracy in estimating community energy costs, efficient congestion management, decreased line flow, improved voltage profile, enhanced reliability and stability, and a reduction in the impact of uncertainties associated with PV generation and demand for energy scheduling.

\section{Conclusion}\label{sec:conclusion}
This study presented a comprehensive energy market that encompasses both virtualization and realization. The day-ahead market focused on virtualization, utilizing a deep learning algorithm to generate quantile scenarios. Subsequently, location scenario-based pricing is introduced to establish a connection, guiding the CO and the DSO in minimizing the impact of forecasting errors related to PV generation and demand. The real-time market addressed the challenges of the realization energy market by considering ESS energy patterns based on actual PV generation and demand data. Payment considerations within the energy market accurately reflected true costs. 

Moreover, the simulation results emphasized the strategic significance of household energy scheduling within the distribution system, enabling prediction of potential P2P energy costs and effective congestion management. Noteworthy, the proposed approach led to a substantial reduction in optimality gap within the community across various seasons. In addition, the charging and discharging patterns of ESSs exhibited variations in response to uncertainty, compared to the case where uncertainty was ignored. This holds considerable implications for decision-making and resource allocation within the energy community, providing confidence in both the CO and the DSO to make well-informed decisions that contribute to effective energy management. In the future, we plan to develop machine learning techniques to reduce further the optimiality gap.

\section{appendix}
\subsection{DLMP parameters} \label{appdix:DLMP_calulcation}
\begin{subequations}
\begin{align}
    C1_{t|ud} = \frac{((f_{it|ud}^P)^2 + (f_{it|ud}^Q)^2)x_{i} + l_{it|ud}f_{it|ud}^Q(r_{i}^2-x_{i}^2)}{((f_{it|ud}^P)^2 + (f_{it|ud}^Q)^2)x_{i} - l_{it|ud}f_{it|ud}^Q(r_{i}^2 + x_{i}^2)} \nonumber\\
    - \frac{2l_{it|ud}f_{it|ud}^{P}r_{i}x_{i}}{((f_{it|ud}^P)^2 + (f_{it|ud}^Q)^2)x_{i} - l_{it|ud}f_{it|ud}^Q(r_{i}^2 + x_{i}^2)}
\end{align}
\begin{align}
    C2_{t|ud} = \frac{((f_{it|ud}^P)^2 + (f_{it|ud}^Q)^2)r_{i} - l_{it|ud}f_{it|ud}^Q(r_{i}^2 + x_{i}^2)}{((f_{it|ud}^P)^2 + (f_{it|ud}^Q)^2)x_{i} - l_{it|ud}f_{it|ud}^Q(r_{i}^2 + x_{i}^2)}
\end{align}
\begin{align}
    C3_{t|ud} = \frac{((f_{it|ud}^P)^2 + (f_{it|ud}^Q)^2)x_{i} + l_{it|ud}f_{it|ud}^Q(r_{i}^2-x_{i}^2)}{((f_{it|ud}^P)^2 + (f_{it|ud}^Q)^2)x_{i} - l_{it|ud}f_{it|ud}^Q(r_{i}^2 + x_{i}^2)} \nonumber\\
    + \frac{2l_{it|ud}f_{it|ud}^{Q}r_{i}x_{i}}{((f_{it|ud}^P)^2 + (f_{it|ud}^Q)^2)x_{i} - l_{it|ud}f_{it|ud}^Q(r_{i}^2 + x_{i}^2)}
\end{align}
\begin{align}
    C4_{t|ud} = \frac{2((f_{it|ud}^P)^{3}r_{i} - (f_{it|ud}^Q)^{3}x_{i})} {((f_{it|ud}^P)^2 + (f_{it|ud}^Q)^2)x_{i} - l_{it|ud}f_{it|ud}^Q(r_{i}^2 + x_{i}^2)} \nonumber\\ 
    + \frac{2f_{it|ud}^{P}f_{it|ud}^Q(f_{it|ud}^{P}r_{i}-f_{it|ud}^{Q}x_{i})}{((f_{it|ud}^P)^2 + (f_{it|ud}^Q)^2)x_{i} - l_{it|ud}f_{it|ud}^Q(r_{i}^2 + x_{i}^2)}
\end{align}
\begin{align}
    C5_{t|ud} = \frac{2((f_{it|ud}^P)^{3}r_{i} - (f_{it|ud}^Q)^{3}x_{i})} {((f_{it|ud}^P)^2 + (f_{it|ud}^Q)^2)x_{i} - l_{it|ud}f_{it|ud}^Q(r_{i}^2 + x_{i}^2)} \nonumber\\ 
    + \frac{2f_{it|ud}^{P}f_{it|ud}^Q(f_{it|ud}^{P}r_{i}-f_{it|ud}^{Q}x_{i})}{((f_{it|ud}^P)^2 + (f_{it|ud}^Q)^2)x_{i} - l_{it|ud}f_{it|ud}^Q(r_{i}^2 + x_{i}^2)} \nonumber \\
    + \frac{2l_{it|ud}^2((f_{it|ud}^P)r_{i}^3 + (f_{it|ud}^Q)x_{i}^3)} {((f_{it|ud}^P)^2 + (f_{it|ud}^Q)^2)x_{i} - l_{it|ud}f_{it|ud}^Q(r_{i}^2 + x_{i}^2)} \nonumber \\
    - \frac{ 4l_{it|ud}^{P}f_{it|ud}^Q(r_{i}^2-x_{i})^2} {((f_{it|ud}^P)^2 + (f_{it|ud}^Q)^2)x_{i} - l_{it|ud}f_{it|ud}^Q(r_{i}^2 + x_{i}^2)} \nonumber\\ 
    + \frac{4l_{it|ud}r_{i}x_{i}((f_{it|ud}^{P})^2-(f_{it|ud}^{Q})^2)}{((f_{it|ud}^P)^2 + (f_{it|ud}^Q)^2)x_{i} - l_{it|ud}f_{it|ud}^Q(r_{i}^2 + x_{i}^2)} \nonumber\\ 
    + \frac{2l_{it|ud}^{2}r_{i}x_{i}((f_{it|ud}^{P})r_{i}-(f_{it|ud}^{Q})x_{i})}{((f_{it|ud}^P)^2 + (f_{it|ud}^Q)^2)x_{i} - l_{it|ud}f_{it|ud}^Q(r_{i}^2 + x_{i}^2)}
\end{align}
\end{subequations}

\renewcommand{\baselinestretch}{0.85}
\bibliographystyle{IEEEbib}

\end{document}